\documentclass[11pt,a4paper, reqno]{amsart}
\usepackage{amssymb, paralist}
\usepackage{amsfonts, amsmath, wasysym}
\usepackage{url}

\usepackage{calligra}
\usepackage{bbm}
\usepackage{mathtools}

\usepackage{latexsym}
\usepackage{graphicx, color}
\usepackage{xcolor}
\usepackage{indentfirst}
\usepackage{stmaryrd}
\usepackage[mathscr]{euscript}

\usepackage{amssymb}
\usepackage{amsmath, amscd}
\usepackage{graphicx, color}
\usepackage{indentfirst}

\newcommand{\N}{\mathbb{N}} 

\newcommand{\R}{\mathbb{R}}

\newcommand{\ds}{\displaystyle}
\renewcommand{\>}{\right>}
\newcommand{\flecha}{\longrightarrow}

\newcommand{\ii}{\mathbbm}

\newcommand{\mc}[1]{\mathcal{#1}}

\newcommand{\E}{\mathbb E}

\def\fle{\rightarrow}
\def\parcial#1#2{\fracc{\partial #1}{\partial#2}}
\def\parci#1#2{\tfrac{\partial #1}{\partial#2}}

\def\deri#1#2{\fracc{d #1}{d#2}}

\def\({\left (}
\def \){\right)}
\newtheorem{mainthm}{Theorem}[]
\newtheorem{maincor}[mainthm]{Corollary}
\newtheorem{thm}{Theorem}[section]

\newtheorem*{thm*}{Theorem} 
\newtheorem{cor}[thm]{Corollary}\newtheorem{prop}[thm]{ Proposition}

\theoremstyle{definition}

\newtheorem{defin}[thm]{Definition}

\newtheorem*{rems*}{Remarks}
\theoremstyle{remark}
\newtheorem{rem}[thm]{Remark}

\newcommand{\eps}{\ensuremath{\varepsilon}}

\newtheorem{teor}{\hspace{12pt} Theorem}

\newtheorem{lema}[teor]{\hspace{12pt} Lemma}

\numberwithin{teor}{section}

\newcommand{\be}{\begin{enumerate}}
\newcommand{\ee}{\end{enumerate}}
\newcommand{\bi}{\begin{itemize}}
\newcommand{\ei}{\end{itemize}}
\newcommand{\bd}{\begin{description}}
\newcommand{\ed}{\end{description}}
\newcommand{\bec}{\begin{equation}}
\newcommand{\eec}{\end{equation}}
\newcommand{\ba}{\begin{array}}
\newcommand{\ea}{\end{array}}
\newcommand{\bt}{\begin{thm}}
\newcommand{\et}{\end{thm}}
\newcommand{\bdem}{\begin{proof}}
\newcommand{\edem}{\end{proof}}
\newcommand{\bl}{\begin{lema}}
\newcommand{\el}{\end{lema}}
\newcommand{\bnp}{\begin{rem}}
\newcommand{\enp}{\end{rem}}
\newcommand{\bde}{\begin{defin}}
\newcommand{\ede}{\end{defin}}
\newcommand{\bnod}{\begin{rem}}
\newcommand{\enod}{\end{rem}}
\newcommand{\bp}{\begin{prop}}
\newcommand{\ep}{\end{prop}}
\newcommand{\bco}{\begin{cor}}
\newcommand{\eco}{\end{cor}}

\newcommand{\nn}{\nonumber}
\newcommand{\lb}{\label}

\newcommand{\G}{\mathsf{G}}

\DeclareMathOperator{\Ric}{Ric}

\newcommand{\ene}{\end{equation} }

\newcommand{\fint}{{-}\hspace*{-1.05em}\int}

\def\G{\Gamma}
\def\b{\beta}

\DeclareRobustCommand{\rchi}{{\mathpalette\irchi\relax}}
\newcommand{\irchi}[2]{\raisebox{\depth}{$#1\chi$}}

\def\a{\alpha}
\def\g{\gamma}

\usepackage{bm}

\usepackage[draft]{optional}

\renewcommand{\(}{\left(}

\renewcommand{\>}{\right>}
\renewcommand{\)}{\right)}
\def\bal{\begin{align}}
\def\eal{\end{align}}

\newcommand{\st}{M \times I}

\numberwithin{equation}{section}
\def\be{\begin{equation}}
\def\ee{\end{equation}}

\def\a{\alpha}

\def\parcial#1#2{\frac{\partial #1}{\partial#2}}

\def\deri#1#2{\frac{d #1}{d#2}}

\def\flecha{\longrightarrow}
\def\fle{\rightarrow}

\def\ds{\displaystyle}

\newcommand{\gi}{G^{\p 0 \p 0}}

\newcommand{\p}{\scriptscriptstyle}

\newcommand{\pa}{\p{\text{par}}}

\begin{document}

\markright{}


\pagestyle{myheadings}

\title{Brownian motion on Perelman's almost Ricci-flat manifold}

\author{Esther Cabezas-Rivas}
\address{Goethe-Universit\"at Frankfurt, Robert-Mayer-Strasse~10, 60325 Frankfurt, Germany}
\email{cabezas-rivas@math.uni-frankfurt.de}

\author{Robert Haslhofer}
\address{Department of Mathematics, University of Toronto, 40 St George Street, Toronto, ON M5S 2E4, Canada
}
\email{roberth@math.toronto.edu}

\thanks{ The first named author has been partially supported by by  the
	MINECO (Spain) and FEDER  project MTM2016-77093-P.
	The second named author   has been partially supported by NSERC grant RGPIN-2016-04331, NSF grant DMS-1406394 and a Connaught New Researcher Award.
	We are indebted to Peter Topping for the proposal of applying the results in \cite{Na13} to Perelman's manifold in order to recover the corresponding statements in \cite{HaNa}, which inspired the present paper.
	We would like to thank {\it The Fields Institute for Research in Mathematical Sciences} for providing an excellent research atmosphere and hosting the Thematic Program on Geometric Analysis from July--December 2017, when this work was completed.}

\maketitle

\begin{abstract}
We study Brownian motion and stochastic parallel transport on Perelman's almost Ricci flat manifold $\mathscr M=M\times \mathbb S^{\p N}\times I$, whose dimension depends on a parameter $N$ unbounded from above. We construct sequences of projected Brownian motions and stochastic parallel transports which for $N \to \infty$ converge to the corresponding objects for the Ricci flow. In order to make precise this process of passing to the limit, we study the martingale problems for the Laplace operator on $\mathscr M$ and for the horizontal Laplacian on the orthonormal frame bundle $\mathscr{OM}$ .

As an application, we see how the characterizations of two-sided bounds on the  Ricci curvature established by A.~Naber applied to Perelman's manifold lead to the inequalities that characterize solutions of the Ricci flow discovered by Naber and the second author.
\end{abstract}

\section{Introduction}

Typically, elliptic theories are regarded as the static case of the corresponding parabolic problem; in that sense, many times the better-understood elliptic theory has been a source of intuition to generalize the corresponding results in the parabolic case. Examples of this feedback are minimal surfaces/mean curvature flow and harmonic maps/solutions of the heat equation. 

In the present paper, we will focus on an  alternative approach: view the parabolic theory as a limiting case of the related elliptic problem as the dimension goes to infinity. This interpretation appears in the previous literature (see \cite[\S 6]{Per1}, \cite[Section 8]{LottOT} and \cite[Section 2.10.2]{Taob}), but usually is only treated from a heuristic point of view.

The approximation of a parabolic problem by taking suitable limits of an elliptic counterpart has been fruitfully used before, for instance, this idea is behind elliptic regularization techniques. In particular, a possible way to construct weak solutions of the mean curvature flow is based on the level set flow method \cite{CGG,ES}, which can be geometrically interpreted as solving an associated elliptic problem in an ambient space that is one dimension higher \cite{Ilmanen}.

To find a similar connection between almost Einstein manifolds and Ricci flow, it is not that clear how to implement analogous ideas, because of the tensor nature of the equations involved and the lack of an ambient space.  The present paper succeeds to provide such a correspondence by performing analysis on the path space of  Perelman's almost Ricci flat manifold $\mathscr M$ (see \cite[\S 6]{Per1}), whose dimension depends on a parameter $N$ unbounded from above. After a suitable averaging and limiting process as $N \to \infty$ we get the corresponding analytic properties for the path space associated to a manifold evolving under the Ricci flow.

 Let us briefly describe Perelman's manifold: Let $(M, g(t))_{t \in [0, T_1)}$ be a solution of the Ricci flow, which once seen with respect to a reverse time parameter $\tau= \mc T - t$ (for some $\mc T \in (0,T_1)$) is defined for $\tau$ lying in some interval $I \subset \mathbb{R}_+$. Perelman considered the manifold 
\begin{equation} \lb{defPerM}
\mathscr{M} = M \times \mathbb S^{\p N} \times  I  \, \text{equipped  with } \,  G = g_\tau  + \tau h+ \Big(\frac{N}{2 \tau} + R\Big) d\tau^2 ,
\end{equation}
where $g_\tau=g(\mc T-\tau)$, $R$ denotes the scalar curvature of $g_\tau$, and $h$ is the metric on $\mathbb S^{\p N}$ with constant curvature $\frac1{2N}$. One can find a constant $C < \infty$ so that
\begin{equation*}
|\Ric_G| \leq \frac{C}{N}.
\end{equation*}
 As these bounds converge to 0 as $N \to \infty$,  we regard $\mathscr M$ as an almost Ricci-flat manifold for very large $N$. This was used by Perelman in \cite{Per1} to heuristically deduce the monotonicity of the reduced volume as a limiting estimate obtained after formally applying the Bishop-Gromov comparison theorem to $(\mathscr M, G)$.

Broadly speaking, our first main result consists in studying the Brownian motion on Perelman's manifold $\mathscr M$ and showing that, after averaging over the sphere factor to remove the $\mathbb S^{\p N}$-dependence and projecting on the space-time part, we obtain a sequence that converges as $N \to \infty$ to the parabolic Brownian motion associated to the family of manifolds evolving under the Ricci flow (Theorem \ref{MainT1}).

Let us illustrate why this result is rather surprising: Even after averaging out the influence of the sphere factor, curves of Brownian motion are allowed to move randomly in any direction of the whole space-time. However, the parabolic Brownian motion associated to a Ricci flow only consists of space-time curves moving backwards in time with unit speed, see \cite{HaNa}. Interestingly enough, our computations reveal that Perelman's metric enforces exactly this type of motion for $N\to \infty$. Namely, on the one hand 
the influence of the big sphere in \eqref{defPerM} causes a large drift term in backwards time direction, and on the other hand the coefficient of Brownian motion in time direction converges to zero for $N\to \infty$. Combining these two effects, for $N\to \infty$ the motion concentrates exactly on those space-time curves that move backwards in time with unit speed.

In addition to Brownian motion, the second main ingredient for doing analysis on path space is stochastic parallel transport. Roughly speaking, the Brownian motion (Wiener measure) allows one to do integration on path space, and the stochastic parallel transport enables one to define suitable derivatives on path space. Our second main result shows the parabolic stochastic parallel transport on a Ricci flow space-time can be obtained as a limit for $N\to \infty$ of a sequence of stochastic parallel transport maps on Perelman's manifold $\mathscr M$ after carrying out an averaging and (somewhat more involved) projection procedure (Theorem \ref{MainT2}).

The machinery from Theorem \ref{MainT1} and Theorem \ref{MainT2} provides a rigorous link between the elliptic and parabolic theories associated to manifolds with bounded Ricci curvature and the Ricci flow, respectively. As an application to illustrate this effective link, we  consider the following two recent developments: the characterizations of two-sided bounds on the Ricci curvature established by A. Naber in \cite{Na13} (see also \cite{ChTh,HaNa_m, WW16}) and the charachterizations of solutions of the Ricci flow obtained by Naber and the second author \cite{HaNa}. By applying the characterizations from \cite{Na13} to Perelman's manifold  $\mathscr M$   and letting $N \to \infty$, we obtain the estimates from \cite{HaNa}, which characterize solutions of the Ricci flow (see Corollary \ref{Cor3}).

\subsection{Background on characterizations of bounded Ricci curvature}

On a smooth Riemannian manifold $M$, it is well known that by exploiting the Bochner formula one can obtain estimates that characterize lower bounds on the Ricci curvature. For instance, if $H_t $ denotes the heat flow, one can easily prove the following characterization via a sharp gradient estimate:
\begin{equation} \lb{grad_est_Ric_fin}
\Ric \geq -\kappa  \quad  \Leftrightarrow \quad |\nabla H_t f| \leq e^{\frac{\kappa t}{2}} H_t |\nabla f|
\end{equation}
for all test functions $f: M \fle \R$. Such characterizations are very useful to make sense of Ricci curvature bounded from below by $-\kappa$ in non-smooth metric measure spaces (cf.~\cite{ AGS,AGS2,EKS,LottV, Sturm}). Apart from the above gradient estimate, one can characterize lower Ricci bounds via optimal transport, by means of sharp logarithmic Sobolev and Poincar\'e inequalities, etc (see \cite{BE85, BL,CCS,RS}).

Until the recent appearance of \cite{Na13} nothing was known about analogous characterizations of two-sided bounds on the Ricci curvature. The main new insight by Naber in \cite{Na13} which allows to characterize also upper bounds on the Ricci curvature consists in understanding the analytic properties of the infinite dimensional path space $P(M)$  instead of working on $M$ itself. Recall that $P(M)$ is the space of continuous curves $X: [0, \infty) \fle M$.
 It is naturally endowed with a family of probability measures $\{\mathbb P_x\}_{x \in M}$, where each $\mathbb P_x $ is the Wiener measure of Brownian motion starting at $x \in M$. 
One can introduce a notion of stochastic parallel transport (see \cite{Elw, Ito, Mall}), which facilitates the definition of a family of gradients $\{\nabla_{s}^{||}\}_{s \geq 0}$ which, loosely speaking,  take derivatives of functions $F: P(M) \fle \R$ along vector fields that are parallel after time $s$.
In \cite{Na13} Naber discovered an infinite-dimensional generalization of \eqref{grad_est_Ric_fin}, which actually characterizes two-sided bounds of the Ricci curvature by means of a sharp gradient estimate:
\begin{multline} \lb{grad_est_Ric_inf}
-\kappa\leq \Ric\leq \kappa \quad \Leftrightarrow  \\
\bigg|\nabla_x \int_{P(M)}  \!\!\! F \, d \mathbb P_x\bigg| \leq \int_{P(M)}\!\!\! \Big(\big|\nabla_{0}^{||}  F\big|  + \frac{\kappa}{2} \int_0^\infty \!\! e^{\kappa s/2} \big|\nabla_{s}^{||}  F\big|\, ds\Big)\,d  \mathbb P_x
\end{multline}
for all test functions $ F: P(M) \fle \R$. To understand the estimate in \eqref{grad_est_Ric_inf}, note that in the special case $F(X) =f(X_t)$, where $f: M \fle \R$ and $t$ are fixed, we recover the finite dimensional estimate \eqref{grad_est_Ric_fin}. Of course there are many more test functions $F$ on path space, and this is one of the reasons why the estimate \eqref{grad_est_Ric_inf} is actually strong enough to characterize two-sided Ricci bounds and not just lower bounds. Naber also proved  infinite dimensional log-Sobolev and Poincar\'e inequalities on path space which also characterize $-\kappa\leq \Ric\leq \kappa$. 

\subsection{Background on characterizations of solutions of Ricci flow} \lb{BackRF}

Let now $(M, g(t))_{t \in [0, T_1)}$ be a family of evolving manifolds. The parabolic counterpart of \eqref{grad_est_Ric_fin}, a characterization of supersolutions of the Ricci flow via sharp gradient estimates, can be obtained by playing around with the Bochner formula: If we write $H_{sT} f$ for the solution at time $T$ of the heat equation $\partial_t w  = \Delta_{g(t)} w$ with initial condition $f$ at time $s$, then it holds that
\begin{equation} \lb{super_RF}
\parci{g}{t}  \geq - 2 \Ric_{g(t)} \quad \Leftrightarrow \quad |\nabla H_{sT} f| \leq H_{sT} |\nabla f|
\end{equation}
for all test functions $f: M \fle \R$.

To prove an infinite dimensional version of \eqref{super_RF}, Naber and the second author introduced in \cite{HaNa} suitable parabolic versions of the notions of path space, Wiener measure, stochastic parallel transport, and parallel gradient. The first main difference with the elliptic case is that on the path space associated to an evolving family of manifolds the curves are not allowed to move freely in all directions of the whole space-time $\mathbb M=M \times [0, T_1]$. Instead, only those curves moving backwards in time with unit speed are taken into account to define the parabolic path space
\begin{equation*}
P^{\pa}_{{\p T}}(\mathbb M)=\big\{(x_s, T-s)\, | \,  s \in [0, T] \text{ and } x \in C^0([0, T], M)\} .
\end{equation*}
This can be endowed with a parabolic Wiener  measure ${\mathbb P}^{\pa}_{\p{(x, T)}}$, given by Brownian motion on $(M, g(t))$ based at $(x, T)\in \mathbb M$. The space-time $\mathbb{M}$ comes equipped with a natural metric connection defined by
\begin{equation} \lb{st_conn}
\nabla_X Y = \nabla_X^{g(t)} Y \qquad \text{and} \qquad \nabla_{t} Y = \partial_t Y + \tfrac1{2}\parci{g}{t}(Y, \, \cdot\,)^{\sharp_{\p g(t)}}.
\end{equation}
Using this space-time connection a notion of stochastic parallel transport on $\mathbb{M}$, and the parabolic parallel gradients $\{^{\pa}\nabla^{||}_s\}_{s \geq 0}$ have been introduced in \cite{HaNa}, and an infinite dimensional generalization of \eqref{super_RF} has been established:
\begin{equation} \lb{grad_est_RF}
\parci{g}{t} = - 2 \Ric_{g(t)} \quad \Leftrightarrow \quad  \bigg|\nabla_x \int_{P^{\pa}_{\p{T}}(\mathbb M)}  F \, d \mathbb P^{\pa}_{\p{(x, T)}}\bigg| \leq \int_{P^{\pa}_{\p T}(\mathbb M)} \big|^{\pa}\nabla_{\p 0}^{||}  F\big| \, d \mathbb P^{\pa}_{\p{(x, T)}}
\end{equation}
for all test functions $F: P_{\p{T}}^{\pa}(\mathbb M) \fle \R$. 
In \cite{HaNa} one can also find corresponding characterizations in terms of log-Sobolev and Poincar\'e inequalities.

\subsection{Main results} \lb{main}

Given a Ricci-flow $(M, g(t))_{t \in [0, T_1)}$, we can associate to it Perelman's almost Ricci-flat manifold $\mathscr M$ as defined in \eqref{defPerM}. With the aim of relating analysis on the path space $P(\mathscr M)$ with analysis on the parabolic path space $P_{\p{T}}^{\pa}(\mathbb M)$, we consider the natural projection $\pi: \mathscr M \fle \st$, and the induced projection map $\hat \pi: P(\mathscr M) \fle P(\st)$ between path spaces. 

Given $p =(x, y, T) \in \mathscr M$, let $\mathbb P_{p}^{\p N}$ be the Wiener measure on  $P(\mathscr M)$ based at $p \in \mathscr M$ and consider the following probability measures defined by averaging out the sphere factor
\begin{equation*}
\overline{\mathbb P}_{\p{(x, T)}}^{\p N} =\fint_{\pi^{-1}(x, T)} \mathbb P_{p}^{\p N} \, d\mu_{_{\pi^{-1}(x, T)}}(p).
\end{equation*}
 After pushing forward the above measure by means of $\hat \pi$, we obtain a sequence
 \bec\lb{eq_seq_for_t1}
{\mathbb Q}_{\p{(x,T)}}^{\p N}=\hat\pi_\ast \overline{\mathbb P}_{\p{(x, T)}}^{\p N}
 \eec
 of probability measures  on  $P(\st)$. Our first main result concerns the limiting and concentration behavior of this sequence of probability measures. Before stating it, let us point out that we have the natural inclusion
 \bec \lb{nat_incl}
 P_{\p{T}}^{\pa}(\mathbb M)\subset P(M\times [0,T]).
 \eec
 and that in the definition of Perelman's manifold in \eqref{defPerM} one can essentially choose $\mc T=T$ and thus $I=[0,T]$.\footnote{To make things precise, one has to sneak in a $\delta>0$ (see Section \ref{sec2} for details on that).}

 \begin{mainthm} {\bf (Convergence of Brownian motion)} \lb{MainT1}
	For every $(x,T)$, the sequence of probability measures $\{{\mathbb Q}_{\p {(x,T)}}^{\p N}\}_{N \in \mathbb N}$ converges to the parabolic Wiener measure $\mathbb P^{\pa}_{\p{(x, T)}}$. In particular, the limiting measure  concentrates on space-time curves that move backwards in time with unit speed.
\end{mainthm}

Having obtained the limiting behavior of Brownian motion for $N\to \infty$, the other main goal will be to show that the stochastic parallel transport on $\mathscr M$ converges in a suitable sense to the corresponding parabolic one on space-time. To this end, we need to recall the construction of the  horizontal Brownian motion $\mc U_s$ on the orthonormal frame bundle $\mathscr{P}: \mathscr{OM} \fle \mathscr{M}$.
Given $\mathscr U \in \mathscr{OM}$, following Eells-Elworthy-Malliavin one solves the 
the stochastic differential equation
\begin{equation} \lb{dU_M}
d \, \mc  U_s = \sum_{{\p A = 0}}^{{\p D-1}} H_{{\p A}}(\mc U_s) \circ dW^{{\p A}}_s \qquad \text{with} \qquad \mc U_0 = \mathscr U. 
\end{equation}
Here $D = n + N + 1$, $H_{\p A}$ are the canonical horizontal fields, $W_s$ is Brownian motion in $\mathbb{R}^D$, and $\circ$ indicates that the equation has to be integrated in the Stratonovich sense.

Now let $\mathbb P_0$ the Wiener measure based at the origin on $P(\R^{\p D})$.
 The solution of \eqref{dU_M} gives a map $\mc U: P_{\p 0} (\R^{\p D}) \fle P_{\mathscr U} (\mathscr{OM})$ that defines a probability measure
\begin{equation} \lb{diff_meas_static}
\mathbb P^{\p N}_{\p{\mathscr U}}= \mc U_\ast \mathbb P_0,\,
\end{equation}
the Wiener measure of horizontal Brownian motion on $\mathscr{OM}$. Note that the Wiener measure of Brownian motion on $\mathscr M$ can be recovered as 
\begin{equation} \lb{rel_P_Pu}
 \mathbb P_p^{\p N} =  \hat{\mathscr P}_{\ast}\mathbb P^{\p N}_{\p{\mathscr U}},
\end{equation}
where $ \hat{\mathscr P}$ denotes the map between path spaces induced by the projection $\mathscr P$.

Let us also recall the parabolic counterpart of the above, which has been constructed in \cite{HaNa}:
Let $\mathfrak p: \mathscr {O} \fle  \mathbb M$  the orthonormal frame bundle  whose fibers are the linear isometries $u: \R^n \fle (T_x M,g_t)$.   On the bundle $\mathscr O$ one solves
\begin{equation} \lb{SDE_RF}
d U_s = -D_s \, ds + \sum_{a = 1}^n H_a(U_s) \circ dW^a_s, \quad \text{with }   U_{\p 0} =  u\in \mathscr{O}_{(x,T)},
\end{equation}
where $(H_1,\ldots,H_n,D_s)$ are the $n+1$ canonical horizontal vector fields for the space-time connection \eqref{st_conn}.
The solution $\{U_s\}_{s \in [0, T]}$ satisfies
$
\frak p(U_s) = (X_s, T-s) \in P^{\pa}_{\p{ T}}(\mathbb M)
$,
see \cite[Proposition 3.7]{HaNa}, and gives a map $U: P_{\p 0}(\R^n)\fle P_{u}(\mathscr O)$, 
which yields a family of probability measures on $P(\mathscr O)$ defined by
\begin{equation} \lb{defQu}
\mathbb P_u^{\pa} =   U_{\ast} \mathbb P_{\p 0},
\end{equation}
the parabolic Wiener measure of horizontal Brownian motion on space-time.
We can recover the  parabolic Wiener measure on $P^{\pa}_{\p{ T}}(\mathbb M)$ via 
\bec \lb{rel_P_Pu_par}
{\mathbb P}^{\pa}_{\p{(x, T)}} = \hat{\frak p}_\ast \mathbb P^{\pa}_u.
\eec 

With the goal of expressing ${\mathbb P}^{\pa}_{\p{(x, T)}}$ as a suitable limit of $\mathbb P^{\p N}_{\p{\mathscr U}}$, we define the following bundles and projection maps: Let $\mathscr F \mapsto M \times I$ be the ${\rm GL}_n$-bundle with fibres given by invertible linear maps $u:\mathbb{R}^n\to T_x M$.
Let $\pi^1: \mathscr M \fle M$ be the projection from Perelman's manifold onto the first factor, and let $\iota: \R^n \fle \R^{{\p D}}$ be the inclusion.
We define the map
\begin{equation*}
\Theta: \mathscr{OM}  \fle  \mathscr  F,\quad \mc U \mapsto  \pi^1_\ast  \circ \mc U \circ \iota ,
\end{equation*}
and denote by
\begin{equation*}
\hat \Theta: P(\mathscr{OM}) \fle P(\mathscr F)
\end{equation*}
 the induced map between the path spaces.
 
Now, for $u \in \mathscr O_{(x,T)}\subset \mathscr F_{(x,T)}$ we define  averaged probability measures  by 
\begin{equation*}
\overline{\mathbb P}_{u}^{\p N} =\fint_{\Theta^{-1}(u)} \mathbb P_{\p {\mc U}}^{\p N} \, d\mu_{_{\Theta^{-1}(u)}}(\mc U).
\end{equation*}
Pushing forward, we obtain a  sequence of probability measures
\begin{equation*}
\mathbb{Q}^{\p N}_u=\hat\Theta_{\ast} \overline{\mathbb P}_{u}^{\p N}
\end{equation*}
on $P(\mathscr F)$. We are now in the position to state our second main result:

\begin{mainthm} \lb{MainT2}{\bf (Convergence of stochastic parallel transport)}
The sequence $\{\mathbb{Q}^{\p N}_u\}_{N \in \N}$  on $P(\mathscr F)$ constructed above converges for $N\to\infty$ to $\mathbb P_u^{\pa}$, the parabolic Wiener measure of horizontal Brownian motion. In particular, the limiting measure concentrates on the orthonormal frame bundle over space-time.
\end{mainthm}

As mentioned at the beginning, Theorem \ref{MainT1} and Theorem \ref{MainT2} provide an effective link between the theory of spaces with bounded Ricci curvature and the analysis of Ricci flow. To illustrate this, let us show how the infinite dimensional gradient estimate for the Ricci flow can be obtained by this limiting process:

\begin{maincor} \lb{Cor3}
If we apply the gradient estimate \eqref{grad_est_Ric_inf} with $\kappa = C/N$ to Perelman's manifold $\mathscr M$ and choose suitable test functions independent of the $\mathbb S^{\p N}$-factor, then after passing to the limit as $N \to \infty$, we obtain the gradient estimate for the Ricci flow:
\begin{equation*}
\bigg|\nabla_x \int_{P^{\pa}_{\p{T}}(\mathbb M)}  F \, d \mathbb P^{\pa}_{\p{(x, T)}}\bigg| \leq \int_{P^{\pa}_{\p T}(\mathbb M)} \big|^{\pa}\nabla_{\p 0}^{||}  F\big| \, d \mathbb P^{\pa}_{\p{(x, T)}}
\end{equation*}
for all $F\in L^2(P^{\pa}_{\p{T}}(\mathbb M))$.
\end{maincor}
\noindent Similarly, we can also obtain the the quadratic variation estimate, the log-Sobolev inequality, and the spectral gap estimate for the Ricci flow by an analogous limiting procedure based on Theorem \ref{MainT1} and Theorem \ref{MainT2} (see Section \ref{Applic}). Further applications of Theorem \ref{MainT1} and Theorem \ref{MainT2} will be discussed elsewhere.

\subsection{Strategy of the proofs}
We hereafter describe the strategy and main steps of the proofs of Theorem \ref{MainT1}, Theorem \ref{MainT2}, and Corollary \ref{Cor3}.

Instead of working directly with the stochastic differential equations describing the problem under consideration, such as \eqref{dU_M} and \eqref{SDE_RF}, we use the martingale approach which was initiated by Stroock and Varadhan \cite{ SVb,SV_1, SV_2}. This formulation in terms of the martingale problem has the advantage of being particularly well suited for compactness and convergence arguments, and is equivalent to solving the related stochastic differential equations in a weak sense.

The starting point to prove Theorem \ref{MainT1} is to regard the Wiener measure $\mathbb P_p^{\p N}$ as the solution of the martingale problem for the Laplacian $\Delta_{\mathscr M}$ based at $p$ (see Section \ref{MGP_Per} for precise definitions). Let us denote by $\mathbb Q_{\p{(x, T)}}$ any subsequential limit of the sequence described in \eqref{eq_seq_for_t1}. We first prove that any such limit (if it exists) is a solution of the martingale problem for the Ricci flow associated to the heat operator $\partial_\tau + \Delta_{g_\tau}$ based at $(x, T)$, see Theorem \ref{lim_sol}. We then use this result and some estimates from its proof to show that the limit indeed exists (Theorem \ref{cptness}), is unique, and concentrates on the parabolic path space, cf.~\eqref{nat_incl}.

To prove Theorem \ref{MainT2} we view the diffusion measure $\mathbb P^{\p N}_{\p{\mathscr U}}$ in \eqref{diff_meas_static} as the solution of the martingale problem associated to the horizontal Laplacian $\Delta_{\mathscr{OM}}$ based at $\mathscr U$ (for preliminaries about computations on the frame bundle see Section \ref{frame_pre}). The basic scheme of proof is the same as for Theorem \ref{MainT1}, but there is a key new difficulty: While the probability measure $\mathbb{P}_u^{\pa}$ from \eqref{defQu} lives on the parabolic path space of the ${\rm O}_n$-bundle $\mathscr{O}\to \mathbb{M}$, the sequence of $\mathbb P^{\p N}_{\p{\mathscr U}}$ lives on a much larger bundle, i.e. depends on many more variables. To facilitate the limiting process, we introduce an intermediate ${\rm GL}_{n+1}$-bundle $\mathscr G \supset \mathscr O$  (see Section \ref{MGP_OPer}). We can pass to a limit on the bundle $\mathscr G$, and the limit satisfies a martingale problem given by some rather complicated differential operator $\mc D + \mathscr N$ (Theorem \ref{lim_sol_fb}). We then, roughly speaking, reduce the number of variables by solving some of the equations. More precisely, plugging some suitable test functions into the martingale problem formula, using that the initial condition is in $\mathscr{O}$, and employing a Gronwall-type argument, we prove that the limiting measure, suitably interpreted, actually concentrates on the ${\rm O}_n$-bundle $\mathscr{O}$. On this smaller bundle, the complicated operator $\mc D+\mathscr N$ simplifies dramatically. Namely, the operator $\mathscr N$ vanishes identically when restricted to $\mathscr O$, and the operator $\mc D$ becomes exactly the horizontal heat operator for the Ricci flow, i.e. $\mc D|_{\mathscr O} = D_\tau + \Delta_{H}$. Using this, it is easy to conclude the proof of Theorem \ref{MainT2}.

Finally, in Section \ref{Applic}, we prove Corollary \ref{Cor3} by applying Nabher's gradient estimate \eqref{grad_est_Ric_inf} with $\kappa = C/N$ to Perelman's manifold $\mathscr M$ and passing to the limit $N\to \infty$ via Theorem \ref{MainT1} and Theorem \ref{MainT2}.


\section{Martingale Problem on Perelman's almost Ricci flat manifold}\label{sec2}

The goal of this section is to prove Theorem \ref{MainT1}. We start with some preliminaries about path space, Perelman's almost Ricci-flat manifold, and the martingale problem (Section \ref{subsec_prel}--\ref{MGP_Per}), and then carry out the proof in Section \ref{subsec_proof1}--\ref{subsec_proof1c}.

\subsection{Convergence on path space} \lb{subsec_prel}

Given a Riemannian manifold $\mc M$, its  path space $P(\mc M)$ is the space of all continuous curves $X:[0,\infty)\to \mc M$ equipped with
the natural topology induced by uniform convergence on bounded intervals. With this topology, $P( \mc M )$ is a Polish space (i.e.~a  separable, completely metrisable, topological space), see \cite[\S 1.3]{SVb}. We denote by $\Sigma=\mathcal{B}_{P(\mc M )}$ the $\sigma$-algebra generated by the Borel sets.

Recall the notion of convergence of
probability measures that we will use repeatedly hereafter: A sequence $\{\mathbb P_k\}_{k \in \N}$ of probability measures on $(P( \mc M),\Sigma)$ converges  to 
a probability measure $\mathbb P$ if 
\bec \lb{conv_pm}
\ds\lim_{k \to \infty}\int \varphi \, d \mathbb P_k = \int \varphi \, d\mathbb P \quad \text{for all bounded continous} \, \varphi: P(\mc M )\to \mathbb{R}.
\eec
There are several characterizations of this convergence 
(sometimes referred to as {\it Portmanteau Theorem}); in particular,  \eqref{conv_pm} is equivalent to 
\bec \lb{Portt}
\limsup_{k \to \infty} \mathbb P_k(\Omega) \leq P(\Omega) \,\,\, \text{for any closed set $\Omega$ in $P(\mc M)$},
\eec
 see e.g. \cite[Theorem 1.1.1]{SVb}.

\subsection{Perelman's almost Ricci flat metric}\label{sec_per_man}

Let $(M,g(t))_{t\in [0,T_1)}$ be a solution of the Ricci flow. To avoid technicalities, we assume throughout the paper that either $M$ is closed or that $(M, g(t))$ is complete for every $t$ and has uniformly bounded curvatures. By standard interior estimates, this also yields bounds for all the derivatives of the curvatures.

Fix $T<T_1$. Let $\delta>0$ small enough so that $\mc T=T+\delta < T_1$, and let $I=[\delta,T]$. Note that $t\in I$ if and only if $\tau=\mc T-t \in I$.
Recall that Perelman's almost Ricci-flat manifold is defined as
\begin{equation*}
\mathscr{M} = M \times \mathbb S^{\p N} \times I
\end{equation*}
equipped with the metric
\begin{equation*}
  G = g_\tau + \tau h + \Big(\frac{N}{2 \tau} + R\Big) d\tau^2,
\end{equation*}
where $g_\tau=g(\mc T-\tau)$, $R$ denotes the scalar curvature of $g_\tau$, and $h$ is the metric on $\mathbb S^{\p N}$ with constant sectional curvature $1/2N$. Throughout the paper, we tacitly assume that $N$ is large enough to ensure that $G$ is positive definite.
 
We use $i, j,k,\ldots $ to denote coordinates indices on the $M$ factor, $\a, \b, \gamma, \ldots $ for those on
the $\mathbb S^{\p N}$ factor, and we denote with index $0$ the coordinate $\tau$ on $I$. 
Moreover,  we use {\small $I, J, K, \ldots$} to denote indices ranging over all parts $0, 1 \leq i \leq n$ and $1 \leq \a \leq N$, and we use doublestroke symbols $\mathbbm{i, j, k}, \ldots$ to denote indices ranging within $\{0, 1, \ldots, n\}$. Usually we denote with capital or calligraphic letters quantities on $\mathscr M$ and with the corresponding lower case or non-calligraphic letters the same quantities on $M$.

With these conventions, the non-vanishing Christoffel symbols $\G_{\p I \p J}^{\p K}$ of the metric $G$ are:
\bec \lb{Cris}
\ba{l} \Gamma_{ij}^k  =   ^g\!\!\Gamma_{ij}^k \medskip \\
\Gamma_{i{\p 0}}^k  =   R_i^k \medskip \\
\Gamma_{{\p 00}}^k  =  - \frac1{2} \nabla^k R \ea \qquad  \ba{l} \Gamma_{ij}^{\p 0}  =  - \gi R_{ij} \medskip \\
\Gamma_{i{\p 0}}^{\p 0}  =   \frac{\gi}{2} \nabla_i R \medskip \\
\Gamma_{{\p 00}}^{\p 0}  =  \frac{\gi}{2} \(\parci{R}{\tau}  + \frac{R}{\tau}\) - \frac1{2 \tau}\ea
\qquad \ba{l} \Gamma_{\a \b}^\g  = \, ^h\Gamma_{\a \b}^\g \medskip \\
\Gamma_{\a{\p 0}}^\g  =  \frac1{2\tau} \delta^\g_\a \medskip \\
\Gamma_{\a \b}^{\p 0}  =  -\frac{\gi}{2}  h_{\a \b}\ea 
\eec

Moreover, there exists a constant $C < \infty$ such that the Ricci curvature of $(\mathscr M, G)$ satisfies the bound
$$|{\rm Ric}_{G}| \leq \frac{C}{N}.$$
\noindent Detailed computations of the Christoffel symbols and the Ricci curvatures of the metric $G$ can be found for instance in \cite{Wei}.

\subsection{The martingale problem on $\bm{\mathscr M}$} \lb{MGP_Per}

Consider Perelman's almost Ricci-flat manifold $(\mathscr M,G)$. As before, let $P(\mathscr M)$ be the path space of all continuous curves $\mathscr{X}: [0,\infty)\to \mathscr M$ equipped with the Borel $\sigma$-algebra $\Sigma=\mathcal{B}_{P(\mathscr M)}$. The  Borel $\sigma$-algebra of path space coincides with the $\sigma$-algebra generated by all the evaluation maps, i.e.
\begin{equation*}
\Sigma = \sigma\{\mathscr E_s \, |\, 0 \leq s < \infty\},
\end{equation*}
where for every $s \in [0, \infty)$ fixed, $\mathscr E_s$ denotes the evaluation map defined by
\begin{equation*}
\mathscr E_s: P(\mathscr M) \flecha \mathscr M \qquad \mathscr{X} \longmapsto  \mathscr{X}_s.
\end{equation*}
Furthermore, there is a natural filtration 
\begin{equation*}
\Sigma_s = \sigma\{\mathscr E_{\p r} \ |\ 0 \leq r \leq s\}\subset \Sigma,
\end{equation*}
which captures the events observable at time $s$. The corresponding filtration and evaluation maps on $M$ are denoted by $\{\sigma_s\}_{s\geq 0}$ and $\eps_s$, respectively.
Note that the projection map
\begin{equation*}
\pi: \mathscr M \fle \st
\end{equation*}
induces a projection map
\begin{equation*}
\hat \pi: P(\mathscr M) \fle P(\st)
\end{equation*}
so that 
\bec \lb{rel_epi}
\eps_s \circ \hat \pi = \pi \circ  \mathscr E_s.
\eec

Given a point $p \in \mathscr M$, the Wiener measure $\mathbb P_{p}^{\p N}$ on $(P(\mathscr M),\Sigma)$  is characterized as the unique solution of the martingale problem for the Laplacian $\Delta_{\mathscr M}$ based at $p$. More precisely, this means that $\mathbb P_{p}^{\p N}$  satisfies the following two conditions:
\bec \lb{mp_pm_1}
\mathbb P_{p}^{\p N} \big(\mathscr{X}_0 =  p\big) = 1,
\eec
and
\bec \lb{mp_pm_2}
\ds \E^{\mathbb P_{p}^{\p N}} \Big[F(\mathscr X_{b}) - \int_{a}^{b}  \Delta_{\mathscr M} F(\mathscr X_s)\, ds \, \Big|\, \Sigma_{a}\Big] = F(\mathscr X_{a})
\eec 
for all smooth functions with compact support $F \in C^\infty_c(\mathscr M)$ and  all  $b > a \geq 0$.\footnote{In this paper, we use the normalization $dW_s\, dW_s=2\, ds$.}

To be precise, we should also say that we put absorbing boundary conditions to take care of the possibility that the motion reaches the boundary $\partial \mathscr M$. Since in the limit the motion will reach $\partial \mathscr M$ exactly at time $T$, for the ease of presentation in the following we don't explicitly indicate the stopping time $T_{\partial \mathscr M}$.

\subsection{The limiting martingale problem}\lb{subsec_proof1}
Recall from the introduction that we consider the sequence of probability measures
\bec \lb{Def_QNx}
\mathbb Q^{\p N}_{\p{(x, T)}}= \hat\pi_\ast \overline{\mathbb P}_{\p{(x, T)}}^{\p N} \quad \text{where} \quad \overline{\mathbb P}_{\p{(x, T)}}^{\p N} =\fint_{\pi^{-1}(x, T)} \mathbb P_{p}^{\p N} \, d\mu_{_{\pi^{-1}(x, T)}}(p).
\eec
In order to prove that $\mathbb Q^{\p N}_{\p{(x, T)}}$ converges to $\mathbb P^{\pa}_{\p{(x, T)}}$ we first identify the martingale problem that any subsequential limit satisfies.

\bt \lb{lim_sol}
Any subsequential limit $\mathbb Q_{\p{(x, T)}}$ of the sequence defined by \eqref{Def_QNx} solves the martingale problem for the Ricci flow, that is, the martingale problem based at $(x, T)$ associated to the heat operator $\partial_\tau + \Delta_{g_\tau}$. More precisely, $\mathbb Q_{\p{(x, T)}}$ satisfies the following two properties:

\begin{enumerate}[(a)]

\item $\mathbb Q_{\p{(x, T)}}\big[X_0 = (x, T)\big] = 1$,

\item $\ds \E^{\mathbb Q_{\p{(x, T)}}}\Big[f(X_{b}) - \int_{a}^{b} (\partial_\tau + \Delta_{g_\tau}) f(X_s)\, ds\, \Big|\, \sigma_{a}\Big] = f(X_{a})$ for all $f \in C^\infty_c(\st)$.
\end{enumerate}
\et

\bdem
If $\mc D$ is a differential operator on a manifold $\mc M$, and $\varphi \in C^\infty_c(\mc M)$, we set
\begin{equation*}
{\mc Z}^{\varphi}_{s}(\mc D):= \varphi(\mathscr X_s) - \int_0^s (\mc D \varphi)(\mathscr X_{r})\, dr \qquad\quad (s\geq 0).
\end{equation*}
Notice that assertion (b) can be rephrased by saying that 
\bec \lb{newb}
\mc Z^f_s\big(\partial_\tau + \Delta_{g_\tau}\big) \text{ is a $\mathbb Q_{\p{(x, T)}}$-martingale for all $f \in C^\infty_c(\st)$}.
\eec
Hereafter we fix an arbitrary function $f \in C^\infty_c(\st)$, and $b>a\geq 0$. We divide the proof of \eqref{newb} into several steps:

\noindent{\bf Step 1.} $\ds\E^{\overline{\mathbb P}_{(x,T)}^{\p N}}\big[{\mc Z}^{{\pi^\ast \!f}}_{b}\!(\Delta_{\mathscr M})\, \rchi_{\hat\pi^{-1}(L)}\big] = \E^{\overline{\mathbb P}_{(x,T)}^{\p N}}\big[{\mc Z}^{{\pi^\ast \!f}}_{a}\!(\Delta_{\mathscr M}) \, \rchi_{\hat\pi^{-1}(L)}\big]$ for all $L \in \sigma_{a}$.

\noindent Notice that equation \eqref{mp_pm_2} says that 
${\mc Z}^{\p F}_{s}(\Delta_{\mathscr M})$  is a $\mathbb P_{\p p}^{\p N}$-martingale for all $F \in C^\infty_c(\mathscr M)$. In particular, for any test function $F$ independent of $\mathbb S^{\p N}$, namely if we choose
 $F = \pi^\ast f$, 
we have that ${\mc Z}^{ {\pi^\ast \!f}}_{s}(\Delta_{\mathscr M})$ is a $\mathbb P_{\p p}^{\p N}$-martingale, i.e.
\begin{equation*}
\E^{\mathbb P_{p}^{\p N}}\big[{\mc Z}^{{\pi^\ast \!f}}_{b}\!(\Delta_{\mathscr M})\, \big|\, \Sigma_{a}\big] =  {\mc Z}^{{\pi^\ast \!f}}_{ a}(\Delta_{\mathscr M}).
\end{equation*}
By the definition of conditional expectation this means that
\bec \lb{cond_exp_1}
\E^{\mathbb P_{p}^{\p N}}\big[{\mc Z}^{{\pi^\ast \!f}}_{b}\!(\Delta_{\mathscr M}) \, \rchi_{\mc L}\big] = \E^{\mathbb P_{p}^{\p N}}\big[ {\mc Z}^{{\pi^\ast \!f}}_{a}\!(\Delta_{\mathscr M}) \, \rchi_{\mc L}\big] \quad \text{for all } \mc L \in  \Sigma_{a}.
\eec
To proceed, recall that
\begin{equation*}
\sigma_a = \sigma\big(\{\eps_{s}^{-1}(U)\ |\ s \leq a, \, U \subset \st \text{ open}\}\big).
\end{equation*}
Now for any $s\leq a$ and $U \subset \st$ open, using the composition rule \eqref{rel_epi} for evaluation maps we see that
\begin{equation*}
\hat \pi^{-1} \big(\eps_s^{-1}(U)\big)  =  \mathscr E_s^{-1}(\pi^{-1}(U)) \in \Sigma_a.
\end{equation*}
Since $\Sigma_a$ is a $\sigma$-algebra, we infer that $\hat \pi^{-1}(L) \in \Sigma_a$ whenever $L \in \sigma_a$.
We can thus take $\mc L = \hat \pi^{-1}(L)$ with $L \in \sigma_{a}$ in \eqref{cond_exp_1}. Finally, as we are taking expectations of random variables that are constant on the $\mathbb S^{\p N}$-factor, we can replace ${\mathbb P}_{p}^{\p N}$ by ${\overline{\mathbb P}_{\p{(x, T)}}^{\p N}}$, and this finishes the proof of Step 1.

\noindent{\bf Step 2.} $\Delta_{\mathscr M}(\pi^\ast f)  = \pi^\ast(\mc D^{\p N}f)$ with $\mc D^{\p N}f = (\partial_\tau + \Delta_{g_\tau}) f + \mathscr E_{\p N}$, where $|\mathscr E_{\p N}| \leq C/N$ for some constant $C < \infty$.

\noindent Using expression \eqref{Cris} for the Christoffel symbols, one can  compute that $\Delta_{\mathscr M}$ has the following form when acting on functions $F$ constant on the $\mathbb S^{\p N}$ factor:
\begin{align*}
\Delta_{\mathscr M} F &= G^{{\p IJ}}(\partial_{\p I}\partial_{\p J}-\Gamma_{{\p IJ}}^{\p K}\partial_{\p K})F\\
&= \Delta_{g_\tau}F - \Big(G^{\alpha \beta} \Gamma_{\alpha \beta}^{\p 0} + G^{ij} \Gamma_{ij}^{\p 0} + G^{\p{00}} \Gamma_{\p{00}}^{\p 0}\Big) \parci{F}{\tau} + \gi \tfrac{\partial^2 F}{\partial  \tau^2} - \gi \Gamma^i_{\p 0 \p 0} \parci{F}{x^i}\nonumber
\\ & = \Delta_{g_\tau}F +\Big(\tfrac{\gi}{2\tau} N + \gi R - \tfrac{(\gi)^2}{2}(\parci{R}{\tau} + \tfrac{R}{\tau}) + \tfrac{\gi}{2 \tau}\Big) \parci{F}{\tau} + \gi \tfrac{\partial^2 F}{\partial \tau^2} + \tfrac{\gi}{2} \nabla^i R\parci{F}{x^i}.\nonumber
\end{align*}
Recalling that $G^{00}=\left(\tfrac{N}{2\tau}+R\right)^{-1}$ we infer that $\Delta_{\mathscr M}(\pi^\ast f)  = \pi^\ast(\mc D^{\p N}f)$, where
\begin{equation*}
\mc D^{\p N}f = (\partial_\tau + \Delta_{g_\tau}) f + \mathscr E_{\p N},
\end{equation*}
and
\begin{align*}
\mathscr E_{\p N} =   \Big(\tfrac1{N+2 \tau R} - \tfrac{(G^{\p{00}})^2}{2}(\parci{R}{\tau}  + \tfrac{R}{\tau})\Big) \parci{f}{\tau} +G^{\p{00}} \Big(\tfrac{\partial^2 f}{\partial \tau^2}   + \tfrac1{2} \nabla^i R \parci{f}{x^i}\Big) .
\end{align*}
Observe that
\begin{equation*}
\max \{ G^{00}, (N+2\tau R)^{-1}\} \leq C/N.
\end{equation*}
Furthermore, recall that the curvatures and their derivatives are bounded. Since $f$ is smooth with compact support, its derivatives are also bounded. We thus conclude that
 $|\mathscr E_{\p N}|\leq {C}/{N}$ for some $C<\infty$, as desired.

\noindent{\bf Step 3.} $\ds \E^{\overline {\mathbb P}_{\p{(x, T)}}^{\p N}}\big[{\mc Z}^{{\pi^\ast \!f}}_{s}\!(\Delta_{\mathscr M}) \,  \rchi_{\hat \pi^{-1}(L)}\big] = \E^{\mathbb Q_{\p{(x, T)}}^{\p N}}[\mc Z^{f}_s(\mc D^{\p N}) \, \rchi_L]$ for all $L \in \sigma_s$.

\noindent The composition rule \eqref{rel_epi} for evaluation maps, the identity from Step 2, and the definition of push-forward measure yield
\begin{align*}
\E^{\overline {\mathbb P}_{\p{(x, T)}}^{\p N}}\big[{\mc Z}^{{\pi^\ast \!f}}_{s}\!(\Delta_{\mathscr M}) \,  \rchi_{\hat \pi^{-1}(L)}\big]  & = \int_{P(\mathscr M)}\!\! \Big[ \pi^\ast f \circ \mathscr E_\tau  - \int_0^s  \pi^\ast(\mc D^{\p N } f)\circ \mathscr E_{\p \xi} \, d \xi \Big]\!\cdot\! (\rchi_{L} \circ \hat \pi) \, d\overline {\mathbb P}_{\p{(x, T)}}^{\p N} \nonumber\\
& = \int_{P(\st)} \Big[f \circ \eps_\tau - \int_0^s  \mc D^{\p N } f \circ \eps_{\p \xi} \, d \xi \Big]\!\cdot\! \rchi_L  \, d\big[\hat \pi_\ast\overline {\mathbb P}_{\p{(x, T)}}^{\p N}\big] \nonumber
\\ & = \E^{\mathbb Q_{\p{(x, T)}}^{\p N}}\big[\mc Z^{f}_s(\mc D^{\p N}) \, \rchi_L\big],
\end{align*}
which is what we wanted to prove.

\noindent{\bf Step 4.} Assertion \eqref{newb} holds true.

\noindent  Take $b > a \geq 0$ and  any $L \in \sigma_{a} \subset \sigma_{b}$. If $\mathbb Q_{\p{(x, T)}}^{\p N}$ converges to $\mathbb Q_{\p{(x, T)}}$, then by 
the definition of convergence of probability measures we have
\begin{align} \lb{Step 5}
 \E^{\mathbb Q_{\p{(x, T)}}}\big[\mc Z^{f}_a(\partial_\tau + \Delta_{g(\tau)}) \, \rchi_L\big]  & = \lim_{N \to \infty}  \E^{\mathbb Q_{\p{(x, T)}}^{\p N}}\big[\mc Z^{f}_a(\partial_\tau + \Delta_{g(\tau)}) \, \rchi_L\big] \nn \\ & = \lim_{N\to \infty} \E^{\mathbb Q_{\p{(x, T)}}^{\p N}}\big[\mc Z^{f}_a(\mc D^{\p N}) \, \rchi_L\big],
 \end{align}
where we used that $\ds \big|\mc Z^{f}_a(\mc D^{\p N}) - Z^{f}_a(\partial_\tau + \Delta_{g_\tau})\big| \leq C/N$ by Step 2. Now using equation \eqref{Step 5} we deduce that
\begin{align} \lb{final_thm1}
 \E^{\mathbb Q_{\p{(x, T)}}}\big[\mc Z^{f}_a(\partial_\tau + \Delta_{g_\tau}) \, \rchi_L\big]  & = \lim_{N \to \infty} \E^{\overline {\mathbb P}_{\p{(x, T)}}^{\p N}}\big[{\mc Z}^{{\pi^\ast \!f}}_{a}\!(\Delta_{\mathscr M}) \,  \rchi_{\hat \pi^{-1}(L)}\big]  \,\,\, \text{\small (by Step 3 for $s = a$)}
\nn \\ & = \lim_{N \to \infty} \E^{\overline {\mathbb P}_{\p{(x, T)}}^{\p N}}\big[{\mc Z}^{{\pi^\ast \!f}}_{b}\!(\Delta_{\mathscr M}) \,  \rchi_{\hat \pi^{-1}(L)}\big] \,\,\, \text{\small (by Step 1)}
\nn \\ & = \E^{\mathbb Q_{\p{(x, T)}}}\big[\mc Z^{f}_b(\partial_\tau + \Delta_{g_\tau}) \, \rchi_L\big],
\end{align}
where the last equality follows by reversing the above chain of equalities (including those from \eqref{Step 5}) but using the corresponding steps for $s = b$ instead of $s = a$. As \eqref{final_thm1} holds for any $L \in \sigma_a$ and arbitrary $f\in C^\infty_c(\st)$, by definition of conditional expectation we reach \eqref{newb}.

Finally, it remains to prove (a). As the set $\{X_0 = (x, T)\}$ is $\sigma_0$-measurable and closed, applying \eqref{Portt} and the definition of push-forward measure we get
\begin{align*}
Q_{\p{(x, T)}}\big[X_0 = (x, T)\big] & \geq \limsup_{\p{N \to \infty}} Q^{\p N}_{\p{(x, T)}}\big[X_0 = (x, T)\big] 
\\ & = \lim_{\p{N \to \infty}} \fint_{\pi^{-1}(x, T)} \mathbb P_{p}^{\p N}\big[\mathscr{X}_0 = p\big] \, d\mu_{\pi^{-1}(x, T)}(p)= 1,
\end{align*}
where we used \eqref{mp_pm_1} in the last step. This finishes the proof of the theorem.
\edem

\subsection{Existence of a subsequential limit}\lb{sec_subseq_lim}

We now need to ensure that a subsequential limit indeed exists. Essentially the proof consists  in adapting to our non-Euclidean setting the arguments in Stroock-Varadhan \cite[Section 1.4]{SVb}. The main change is to check that we can construct a large enough family of $\mathbb Q^{\p N}_{\p{(x, T)}}$-submartingales which satisfy the hypotheses of \cite[Theorem 1.4.6]{SVb}. 
\bt \lb{cptness}
 For every $(x, T) \in \st$ there exists a subsequential limit $\mathbb Q_{\p{(x, T)}}$ of the sequence of probability measures  $\big\{\mathbb Q^{\p N}_{\p{(x, T)}}\big\}_{N \in \mathbb N}$ defined by \eqref{Def_QNx}.
\et

\bdem
Given $\rho > 0$ and $S<\infty$, by the precompactness result on path space \cite[Theorem 1.3.2]{SVb}, it is enough to check that
\begin{equation*}
\lim_{\delta \searrow 0} \ds \limsup_{N \to \infty} \mathbb Q^{\p N}_{\p{(x, T)}} \bigg[\sup_{{\tiny \ba{c}0 \leq s_1\leq s_2 \leq S \\ s_2 - s_1 < \delta\ea}} d_{\hat g}(X_{s_1}, X_{s_2}) \leq \rho \bigg] = 1,
\end{equation*}
where we use the metric $\hat g = g_{\tau} + d\tau^2$ on $\st$.

By combining Step 1 and Step 3 from the proof of Theorem \ref{lim_sol}, we deduce that $\mc Z^f_s(\mc D^{\p N})$ is a $\mathbb Q^{\p N}_{\p{(x, T)}}$-martingale for every $f\in C^\infty_c(\st)$, i.e.
\[f(X_a)  =  \E^{\mathbb Q^{\p N}_{\p{(x, T)}}}\Big[f(X_b) - \int_{a}^{b} \mc D^{\p N}  \!f(X_s) \, d s \ \Big| \ \sigma_{a}\Big]\]
for all $b > a \geq 0$. Now by Step 2 from the proof of Theorem \ref{lim_sol} we see that
\bec \lb{def_Af}
\big|\mc D^{\p N}\! f(X_s)\big| \leq \sup_{\st} \big|(\partial_\tau + \Delta_M) f(X_s)\big| + C \leq A_f,
\eec
where $A_f <\infty$ depends only on the curvature bounds for our Ricci flow and the $C^2$-norm of $f$. Accordingly, for any non-negative function $f$ we get
\[f(X_a)  \leq  \E^{\mathbb Q^{\p N}_{\p{(x, T)}}}\Big[f(X_b) + A_f (b - a) \, \Big| \, \sigma_{a}\Big],\]
which says that $(f(X_s) + A_f s, \sigma_s, \mathbb Q^{\p N}_{\p{(x, T)}})$ is a non-negative submartingale. This tells us that $\mathbb Q^{\p N}_{\p{(x, T)}}$ satisfies the first hypothesis in \cite[Theorem 1.4.6]{SVb}. 

The second hypothesis in \cite[Theorem 1.4.6]{SVb} has to be adapted to our non-Euclidean space $\st$. To this end,  instead of  translates of $f$, we use as test functions the  following space-time bump functions: fix $(\bar x, \bar \tau) \in \st$ and choose $f_\rho \in C^\infty_c(\st)$ so that $0 \leq f_\rho \leq 1$, $f_\rho = 1$ on $B_{\frac{\rho}{8}}(p) \times [\bar \tau - \frac{\rho}{8}, \bar \tau + \frac{\rho}{8}]$ and $f_\rho(x, \tau) = 0$ for each $x \in M \setminus B_{\frac{\rho}{4}}(p)$ and each $\tau$ with $|\tau - \bar \tau| \geq \frac{\rho}{4}$. Note that the corresponding constant $A_{f_\rho}$ from \eqref{def_Af} can be chosen independent of the base point $(\bar x, \bar \tau)$. This independence is the suitable replacement for the second hypothesis of \cite[Theorem 1.4.6]{SVb}. The rest of the proof goes through exactly as in \cite[Section 1.4]{SVb}.\edem

\subsection{Concentration on the parabolic path space}\lb{subsec_proof1c}

In order to finish the proof of Theorem \ref{MainT1} it remains to show that our limiting measure concentrates on space-time curves moving backwards in time with unit speed.

\bdem[{Proof of Theorem \ref{MainT1}}] Given the Ricci flow $(M,g(t))_{t\in [0,T_1)}$ and a space-time point $(x,T)\in M\times (0,T_1)$, consider Perelman's manifold $\mathcal{M}=M\times I\times \mathbb S^{\p N}$ equipped with the almost Ricci-flat metric $G$ as in Section \ref{sec_per_man}. By Theorem \ref{lim_sol} and Theorem \ref{cptness} the probability measures $\mathbb{Q}_{{\p (x,T)}}^{\p N}$, which have been defined in \eqref{Def_QNx}, converge along a subsequence to a probability measure $\mathbb{Q}_{{\p (x,T)}}$ which solves the martingale problem based at $(x,T)$ associated to the heat operator $\partial_\tau +\Delta_{g_\tau}$. 

Let $\pi_I: M \times I \fle I$ be the natural projection $\pi_I(x, t) = t$. We apply the conclusion of Theorem \ref{lim_sol} to the particular function $f = h \circ \pi_I$ with $h \in C^\infty_c(I)$ and the curve $X_s =(x_s, \alpha_s)$ with $\alpha \in P(I)$ satisfying $\alpha_0 = T$. Noticing that $(\partial_\tau + \Delta_{g_\tau}) f = - h'$ we get that
\bec \lb{R_mart}
\E^{\mathbb Q_{\p (x, T)}}\Big[h(\alpha_{b}) + \int_{a}^{b}  h'(\alpha_s) \, ds \, \big| \, \sigma_{a}\Big] = h(\alpha_{a})
\eec
for all $b > a \geq 0$ and all $h \in C^\infty_c(I)$.

In particular, choosing $h(t) = t$, $b = s$ and $a = 0$ equation \eqref{R_mart}  yields
\bec\label{eq_conc1}
\E^{{\mathbb Q_{\p (x, T)}}}[\alpha_s] = T -s.
\eec
Next we apply \eqref{R_mart} for the function $h(t) = t^2$, $b = s$, and $a = 0$ to obtain
\begin{align*}
T^2 =\E^{{\mathbb Q_{\p (x, T)}}}\Big[\alpha_{s}^2 + \int_0^s 2 \alpha_r \, dr\Big] & = \E^{{\mathbb Q_{\p (x, T)}}}\big[\alpha_{s}^2\big] + 2 \int_0^s \E^{{\mathbb Q_{\p (x, T)}}}\big[\alpha_r] \, dr \\ & = \E^{{\mathbb Q_{\p (x, T)}}}\big[\alpha_{s}^2 ] - (T -s)^2 + T^2,
\end{align*}
where we used equation \eqref{eq_conc1}. Putting things together, we infer that
\begin{equation*}
\E^{{\mathbb Q_{\p (x, T)}}}\big[\big(\alpha_s - (T-s)\big)^2\big] = \E^{{\mathbb Q_{\p (x, T)}}}\big[\alpha_s^2 - 2 (T-s) \alpha_s + (T-s)^2\big] = 0,
\end{equation*}
and thus $\alpha_s=T-s$ holds $\, \mathbb Q_{\p (x, T)}$-almost surely.

Summing up, $\mathbb Q_{\p (x, T)}$ is a probability measure on $P(M\times [\delta,T])$ supported on space-time curves of the form $\{X_s=(x_s, T-s)\}_{s \in [0, T - \delta]}$ which solves the martingale problem (that we now express using $t$ instead of $\tau$) given by
\begin{align*}
&\mathbb Q_{\p{(x, T)}}\big[X_0 = (x, T)\big] = 1,\\
&\ds \E^{\mathbb Q_{\p{(x, T)}}}\Big[f(x_{b}, T-b) - \int_{a}^{b} \Big[(-\partial_t + \Delta_{g(T  - s)}) f(\cdot, T-s)\Big](x_s)\, ds\, \Big|\, \sigma_{a}\Big] = f(x_{a}, T-a).\nonumber
\end{align*}
Letting $\delta \searrow 0$, this is exactly the martingale problem for the Ricci flow, which is obtained after projecting to $\mathbb M$ the stochastic differential equation \eqref{SDE_RF}. By uniqueness of the solutions of the martingale problem, any subsequential limit ${\mathbb Q_{\p (x, T)}}$ coincides with the parabolic Wiener measure $\mathbb P^{\pa}_{(x, T)}$ and therefore we obtain full convergence (without passing to a subsequence) in Theorem \ref{cptness}.
\edem

\section{Stochastic parallel transport on Perelman's manifold}

The goal of this section is to prove Theorem \ref{MainT2}.

\subsection{Frame bundles over $\bm{\mathscr M}$} \lb{frame_pre}
Set {\small $D= n + N + 1$} and consider the ${\rm GL}_{\p D}$-frame bundle ${\mathscr{P}}: \mathscr{F} \mathscr{M} \fle \mathscr{M}$  whose fibers $ \mathscr{P}^{-1}(p)$ are linear isomorphisms $\mathscr U: \R^{\p D} \fle T_p \mathscr{M}$. The group ${\rm GL}_{\p D}$ acts on ${\mathscr F}\mathscr{M}$ from the right via composition.

Given local coordinates $(x^i, y^\a, \tau)$ on $\mathscr M$,   a local chart on ${\mathscr F} \mathscr M$ is induced as follows. Let $ \parci{}{x^{\p I}} \in \{\parci{}{x^i}, \parci{}{y^\a}, \parci{}{\tau}\}$ be the coordinate frame and denote by $\{E_{\p I}\}_{\p I = \p 0}^{\p{D -1}}$ the canonical basis of $\R^{\p D}$. For a frame $\mathscr U \in {\mathscr F} \mathscr M$ we write
\bec \lb{uea}
\mathscr U  E_{\p A} =  E_{\p A}^{\p I}  \parci{}{x^{\p I}} =  E^i_{\p A} \parci{}{x^i} +  E^\a_{\p A} \parci{}{x^\a} + E^{\p 0}_{\p A} \parci{}{\tau}
\eec
for some matrix $(E_{\p A}^{\p I}) \in {\rm GL}_{\p D}$. Then $(x^i, y^\a, \tau,   E_{\p A}^{\p I}) \in \R^{\p{D + D^2}}$ are local coordinates for ${\mathscr F} \mathscr M$. Note that we are using $A, B, \ldots$ for frame bundle indices, and $I, J, \ldots$ for manifold indices.

The Levi-Civita connection gives a splitting of the tangent bundle
\begin{equation*}
T_{\mathscr U} ({\mathscr F \mathscr M}) = \mc H_{\mathscr U} \oplus \mc V_{\mathscr U}
\end{equation*}
into horizontal and vertical subspaces. The horizontal lift of  $X\in T_{{\mathscr P}(\mathscr U)} \mathscr  M$ to $T_{\mathscr U} ({\mathscr F \mathscr M})$  is the unique vector $X^\ast \in \mc H_{\mathscr U}$ so that ${\mathscr P}_\ast(X^\ast)
= X$.

The orthonormal frame bundle
\begin{equation*}
\mc{O} \mathscr{M}=\big\{\mathscr U \in \mathscr F \mathscr M \ | \ G_{\p{\mathscr P(\mathscr U)}}(\mathscr U E_{\p A}, \mathscr U E_{\p B}) = \delta_{\p A \p B} \text{ for } \text{\small $A, B =0, \ldots, D-1$}\big\}
\end{equation*}
is a subbundle of $\mathscr F \mathscr M$ with  structure group ${\rm O}_{\p D}$. In addition to the canonical horizontal vector fields
\begin{equation*}
H_{\p A}(\mathscr U) = (\mathscr U  E_{\p A})^\ast,
\end{equation*}
there are $D(D - 1)/2$ vertical vector fields defined by
\begin{equation*}
V_{\p A \p B}(\mathscr U) =\deri{}{t}\bigg|_{t = 0} \mathscr U e^{t M_{\p A \p B}},
\end{equation*}
where $M_{\p A \p B} \in \mathfrak{o}(n)$ is the matrix 
 with $-1$ at the {\small $(A, B)$} position, $1$ at the
 {\small $(B, A)$} entry, and 0 elsewhere. We can endow $\mathscr{OM}$ with a Riemannian metric uniquely determined by the formulas
 \begin{equation*}
 \langle X^\ast, Y^\ast \rangle = G(X,Y),\,\,\,\, \langle X^\ast, V_{AB}\rangle = 0,\,\,\,\, \langle V_{AB},V_{CD}\rangle=\delta_{AC}\delta_{BD}-\delta_{AD}\delta_{BC}.
 \end{equation*}
Explicit expressions in coordinates for the horizontal and vertical vector fields are given by
\bec \lb{H_loc}
H_{\p A} =  E_{\p A}^{\p I} \Big( \parci{}{x^{\p I}} - \G_{{\p I} {\p J}}^{\p K} \,  E^{\p J}_{\p B} \,  \parci{}{E^{\p K}_{\p B}}\Big),
\eec
and
\bec \lb{V_loc}
V_{\p{A B}} = E_{\p B}^{\p I} \parci{}{E_{\p A}^{\p I}} - E_{\p A}^{\p I} \parci{}{E_{\p B}^{\p I}},
\eec
respectively, see e.g. \cite[Proposition 2.1.3]{Hsu}.
The {\it horizontal Laplacian} on $\mathscr O \mathscr M$ is defined by 
\bec \lb{def_hor_lap}
\Delta_{_{\mathscr{OM}}} =  \sum_{{\p A} = 0}^{{\p {D -1} }}  H_{\p A}^2.
\eec 
Finally, note that by the definition of  $\mathscr O \mathscr{M}$ and equation \eqref{uea} we have
\begin{equation*}
\delta_{\p A \p B} = \<E_{\p A}^{\p I} \parcial{}{x^{\p I}} , E_{\p B}^{\p J} \parcial{}{x^{\p J}}\>_G = E_{\p A}^{\p I} G_{\p I \p J} E_{\p B}^{\p J}.
\end{equation*}
Equivalently, taking inverses this gives the useful formula
\begin{equation} \label{ginv_ee}
G^{\p I \p J}= \sum_{\p{A = 0}}^{\p{D - 1}} E^{\p I}_{\p A} \, E^{\p J}_{\p A}.
\end{equation}

\subsection{Martingale Problem on the frame bundle of Perelman's manifold} \lb{MGP_OPer}

Consider the path space $P\big(\mathscr O \mathscr M\big)$ with the filtration $\{\Sigma_t\}_{t \geq 0}$ generated by the evaluation maps $\mathscr E_t: P\big(\mathscr O \mathscr M\big) \fle \mathscr O \mathscr M$.  A probability measure $\mathbb P_{\p{\mathscr U}}^{\p N}$ on $P\big(\mathscr O \mathscr M\big)$ is a solution of the martingale problem for the horizontal Laplacian \eqref{def_hor_lap} on $\mathscr O \mathscr M$ based at $\mathscr U \in  \mathscr O \mathscr M$ if the following two conditions hold:
\bec\label{mgp_f1}
\mathbb P_{\p \mathscr U}^{\p N} \big[\,\mc U_0 = \mathscr U\big] = 1,
\eec
and
\bec\label{mgp_f2}
\E^{\mathbb P_{\mathscr U}^{\p N}} \Big[\Psi(\mc U_{s_2}) - \int_{s_1}^{s_2}  \Delta_{_{\mathscr O \mathscr M}} \Psi(\mc U_s)\, ds \, \Big|\,  \Sigma_{s_1}\Big] = \Psi(\mc U_{s_1})\ds
\eec
for all smooth functions with compact support $\Psi \in C^\infty_c(\mathscr O \mathscr M)$ and all $s_2 > s_1 \geq 0$. 

Consider the ${\rm GL}_{n+1}$-bundle $\mathscr{G}$ over $\st$ whose fibers are given by
\bec\lb{def_gln1}
\mathscr{G}_{(x, \tau)}=\{u: \R^{n +1} \flecha \R \times T_x M \, \text{ linear isomorphism}\}.
\eec
Using the inclusion $\iota: \R^{n + 1} \fle \R^D$ and the projection $\varpi: T_{(x, y, T)} \mathscr M \fle T_x M \times \R$, we define the map
\begin{equation*}
\Phi:\mathscr O \mathscr M \flecha \mathscr{G}, \qquad  \mathscr U \longmapsto \varpi \circ \mathscr U \circ \iota,
\end{equation*} 
and we denote by
$\hat \Phi: P\big(\mathscr O \mathscr M\big) \flecha P\big(\mathscr{G}\big)$
the induced map between path spaces.

For every $u \in  \mathscr{O} \subset  \mathscr{G}$, where the inclusion is given by $u \stackrel{\mathfrak i}{\longmapsto} \check u=\left[\begin{array}{c|c}
1 & 0   \\ \hline  0 & u
\end{array}\right]$, we define
\bec \label{Def_QNu}
\mathbb Q^{\p N}_u= \hat\Phi_\ast \overline{\mathbb P}_{u}^{\p N}, \qquad \text{with } \quad\overline{\mathbb P}_{u}^{\p N} =\fint_{\Phi^{-1}(\check u)} \mathbb P_{\p \mathscr U}^{\p N} \, d\mu_{_{\Phi^{-1}(\check u)}}(\mathscr U).
\eec

We now follow the same scheme as in the previous section, and first identify the martingale problem for subsequential limits of the sequence $\{\mathbb Q^{\p N}_u\}_{N \in \mathbb N}$.

\bt \lb{lim_sol_fb}
Any subsequential limit $\mathbb Q_u$ of the sequence defined by \eqref{Def_QNu} solves the martingale problem associated to a certain differential operator $\mathcal D + \mathscr N$ based at $u \in \mathscr O$;  more precisely, $\mathbb Q_u$ satisfies the following two properties:

{\rm (a)} $\mathbb Q_u\big[U_{\p 0} = \check u\big] = 1$,

{\rm (b)} $\ds \E^{\mathbb Q_u}\Big[\psi(U_{s_2}) - \int_{s_1}^{s_2} \big[(\mc D + \mathscr N) \psi\big](U_s)\, ds\, \Big|\, \sigma_{s_1}\Big] = \psi(U_{s_1})$ for all $\psi \in C^\infty_c(\mathscr{G})$,
where $\mc D$ and $\mathscr N$ denote the differential operators
\begin{align} \lb{def_D}
\mc D \psi& = g^{i\ell} \bigg[\tfrac{\partial^2 \psi}{\partial x^\ell \partial x^i} - e^j_b \Big[\parci{\Gamma_{i j}^k}{x^\ell}  \parci{\psi}{e^k_b} + \Gamma_{ij}^k \tfrac{\partial^2 \psi}{\partial x^\ell \partial e^k_b} - e^m_c \G_{\ell m}^r \G_{ij}^k \tfrac{\partial^2 \psi }{\partial e^r_c \partial e^k_b}\Big] + e_c^m \G_{\ell m}^r \G_{ir}^k \parci{\psi}{e^k_c}\bigg]  \nn
\\ & \quad - g^{\ell r} \G_{\ell r}^i \Big( \parci{\psi}{x^i} - \G_{{i} {j}}^{k} \,  e^{j}_{b} \,  \parci{\psi}{e^{k}_{b}}\Big) + \parci{\psi}{\tau}  - e^j_{b} R^k_j \parci{\psi}{e^k_b},
\end{align}
and
\begin{align} \lb{N_ex}
\mathscr N \psi & = \tfrac1{2\tau}  e^{\p 0}_{\ii b} \parci{\psi}{e_{\ii b}^{\p 0}}  +  e_{\p 0}^{j} B^k_j \parci{\psi}{e_{0}^{k}}
 + e_{\ii b}^{\p 0} B^k\parci{\psi}{e_{\ii b}^{k}} \nn
\\ &  \quad - e_{\p 0}^{j} g^{i\ell}\G_{i  j}^{ k}  \Big[\parci{}{x^\ell} - \G_{\ell \ii m}^{r} e_{\ii c}^{\ii m}  \parci{}{e_{\ii c}^r}\Big]\parci{\psi}{e_{0}^{k}} - g^{i\ell} R_{i }^{ k} e^{\p 0}_{\ii b} \Big[\parci{}{x^\ell} - \G_{\ell \ii m}^{r} e_{\ii c}^{\ii m}  \parci{}{e_{\ii c}^r}\Big]\parci{\psi}{e_{\ii b}^{k}}.
\end{align}
Here, $B^k$ and $B^k_j$ denote some bounded functions on $\st$, whose precise form is not needed.
\et

\bdem
Regarding the scheme of proof, we can mimic the steps in the proof of Theorem \ref{lim_sol}, but now working with the martingale problem on the frame bundle given by \eqref{mgp_f1} and \eqref{mgp_f2}. The main difference in carrying this out is that we need to find an explicit expression in this case for the operator $\mc D^{\p N}$ in Step 2. Therefore the goal is to compute $ \Delta_{_{\mathscr{OM}}} (\Phi^\ast \psi)$ for any $\psi \in C^\infty_c(\mathscr G)$. 

Set $\Psi= \Phi^\ast \psi$. This function only depends on the coordinates $x^i$, $\tau$, $E^{\mathbbm k}_{\mathbbm b}$ of $\mathscr{O M}$. Taking this into account we get by means of \eqref{H_loc} that
\begin{equation*}
H_{\p A} \Psi  =  E^{\ii i}_{\p A} \mathscr D_{\ii i} \Psi \qquad \text{with} \qquad \mathscr D_{\ii i}\Psi = \parci{\Psi}{x^{\ii i}} - \G_{\ii i \p J}^{\ii k} E_{\ii b}^{\p J} \parci{\Psi}{E^{\ii k}_{\ii b}}.
\end{equation*}
We can thus write the horizontal Laplacian of $\Psi$ as
\begin{equation*}
\Delta_{_{\mathscr{OM}}} \Psi = \sum_{\p{A = 0}}^{\p{D -1}} H_{\p A}(E^{\ii i}_{\p A} \mathscr D_{\ii i} \Psi) 
= \sum_{\p{A = 0}}^{\p{D -1}}  \Big(E^{\ii i}_{\p A} E^{\ii l}_{\p A}  \mathscr D_{\ii l}\mathscr D_{\ii i} \Psi + H_{\p A}(E^{\ii i}_{\p A}) \mathscr D_{\ii i} \Psi\Big).
\end{equation*}
To proceed, using again \eqref{H_loc} and the fact that $\parci{E^{\ii i}_{\p A}}{x^{\ii j}}   = 0$, we compute
\begin{equation} \label{Haea}
\sum_{\p{A = 0}}^{\p{D -1}}  H_{\p A}(E_{\p A}^{\ii i}) = - \sum_{\p{A = 0}}^{\p{D -1}}   E_{\p A}^{\p J} \,  E_{\p A}^{\p K} \, \Gamma_{\p J \p K}^{\ii i} = - G^{\p J \p K} \Gamma_{\p J \p K}^{\ii i},
\end{equation}
where the last equality follows from \eqref{ginv_ee}. In particular, if we denote by  $\mathscr E_{\p N}$ any general function satisfying $\lim_{N \to \infty} \mathscr E_{\p N} = 0$, and  use the formula \eqref{Cris} for the Christoffel symbols, we get
\begin{equation*}
\sum_{\p{A = 0}}^{\p{D -1}} H_{\p A}(E^{\p 0}_{\p A}) = -(G^{\p 0 \p 0} \Gamma_{\p 0 \p 0}^{\p 0}+  g^{ms} \Gamma_{ms}^{\p 0} + G^{\a\b} \Gamma_{\a\b}^{\p 0}) = 1 + \mathscr E_{\p N},
\end{equation*}
and
\begin{equation*}
\sum_{\p{A = 0}}^{\p{D -1}} H_{\p A}(E^{i}_{\p A}) = - (g^{\ell m}  \Gamma_{\ell m}^i + G^{\p 0 \p 0} \Gamma_{\p 0 \p 0}^i ) = - g^{\ell m}  \Gamma_{\ell m}^i + \mathscr E_{\p N}.
\end{equation*}
Putting things together, the formula for the horizontal Laplacian of $\Psi$ becomes
\begin{align*}
\Delta_{_{\mathscr{OM}}} \Psi  & =  G^{\ii i \ii l} \mathscr D_{\ii l} \mathscr D_{\ii i} \Psi + \mathscr D_0 \Psi - g^{\ell m} \G_{\ell m}^{ i}   \mathscr D_{ i} \Psi + \mathscr E_{\p N}\nonumber
\\ & = g^{i\ell} \Big(\parci{}{x^\ell} - \G_{\ell \ii m}^{r} E_{\ii c}^{\ii m} \parci{}{E_{\ii c}^{r}}\Big)\Big(\parci{\Psi}{x^i} - \G_{i \ii j}^{ k} E_{\ii b}^{\ii j} \parci{\Psi}{E_{\ii b}^{ k}}\Big) \nonumber
\\ & \quad + \parci{\Psi}{\tau} - \G_{{\p 0} {\p J}}^{\ii k} E_{\ii b}^{\p J}\parci{\Psi}{E_{\ii b}^{\ii k}} - g^{\ell m} \G_{\ell m}^i \Big(\parci{\Psi}{x^i} - \G_{i {\p J}}^{\ii k} E_{\ii b}^{\p J}\parci{\Psi}{E_{\ii b}^{\ii k}}\Big) + \mathscr E_{\p N}\nonumber
\\ & = g^{i\ell} \bigg(\tfrac{\partial \Psi}{\partial x^\ell \partial x^i} - E_{\ii b}^{\ii j}\Big[\parci{\G_{i \ii j}^{ k}}{x^\ell}   + \G_{i \ii j}^{ k}  \Big(\parci{}{x^\ell} - \G_{\ell \ii m}^{r} E_{\ii c}^{\ii m}  \parci{}{E_{\ii c}^r}\Big)\Big]\parci{\Psi}{E_{\ii b}^{k}} +  \G_{\ell \ii m}^{r} E_{\ii c}^{\ii m} \G_{i r}^{ k}  \parci{\Psi}{E_{\ii c}^k}\bigg)  \nonumber
\\ & \quad + \parci{\Psi}{\tau} - R^k_j E^j_{\ii b} \parci{\Psi}{E_{\ii b}^k} - \G_{\p{00}}^{\ii k} E^{\p 0}_{\ii b} \parci{\Psi}{E_{\ii b}^{\ii k}} - \G_{\p{0} j}^{\p 0} E^{j}_{\ii b} \parci{\Psi}{E_{\ii b}^{\p 0}} \nonumber
\\ & \quad - g^{\ell r} \G_{\ell r}^i\Big(\parci{\Psi}{x^i} - \G_{i j}^{k} E^{j}_{\ii b} \parci{\Psi}{E_{\ii b}^{k}} - \G_{i \p{0}}^{\ii k} E^{\p 0}_{\ii b} \parci{\Psi}{E_{\ii b}^{\ii k}} - \G_{i j}^{\p 0} E^{j}_{\ii b} \parci{\Psi}{E_{\ii b}^{\p 0}}\Big) + \mathscr E_{\p N}.
\end{align*}
Next we use again \eqref{Cris} and group together  in a new differential operator $\mathscr N$ all the terms which contain an index 0 somewhere. This yields
\begin{equation*}
\Delta_{_{\mathscr{OM}}}\Psi = (\mc D + \mathscr N)\Psi + \mathscr E_{\p N},
\end{equation*}
with
\begin{align*}
 \mc D\Psi  & =  g^{i\ell} \bigg(\tfrac{\partial \Psi}{\partial x^\ell \partial x^i} - E_{b}^{j}\Big[\parci{\G_{i j}^{ k}}{x^\ell}   + \G_{i  j}^{ k}  \Big(\parci{}{x^\ell} - \G_{\ell m}^{r} E_{ c}^{ m}  \parci{}{E_{ c}^r}\Big)\Big]\parci{\Psi}{E_{b}^{k}} +  \G_{\ell  m}^{r} E_{c}^{m} \G_{i r}^{ k}  \parci{\Psi}{E_{c}^k}\bigg) \nonumber
 \\ & \quad + \parci{\Psi}{\tau} - R^k_j E^j_{b} \parci{\Psi}{E_{b}^k} - g^{\ell r} \G_{\ell r}^i\Big(\parci{\Psi}{x^i} - \G_{i j}^{k} E^{j}_{b} \parci{\Psi}{E_{b}^{k}}\Big),
\end{align*}
and
\begin{align*}
\mathscr N \Psi & = \tfrac1{2\tau}  E^{\p 0}_{\ii b} \parci{\Psi}{E_{\ii b}^{\p 0}}  +  E_{\p 0}^{j} \Big[g^{i\ell}
\Big(\G_{i \ell}^r \G_{r j}^{k} + \G_{\ell j}^{r}  \G_{i r}^{ k} -\parci{\G_{i j}^{ k}}{x^\ell} \Big) - R_j^k\Big] \parci{\Psi}{E_{0}^{k}}
\\ & \quad  + E_{\ii b}^{\p 0}\Big[g^{i\ell}
\Big( \G_{i\ell}^r R_{r}^{k}  + R^r_\ell \G_{ir}^k -\parci{R_{i}^{ k}}{x^\ell}  \Big) - \G_{\p{00}}^k\Big]\parci{\Psi}{E_{\ii b}^{k}}
\nonumber\\ &  \quad - E_{\p 0}^{j} g^{i\ell}\G_{i  j}^{ k}  \Big[\parci{}{x^\ell} - \G_{\ell \ii m}^{r} E_{\ii c}^{\ii m}  \parci{}{E_{\ii c}^r}\Big]\parci{\Psi}{E_{0}^{k}} - g^{i\ell} R_{i }^{ k} E^{\p 0}_{\ii b} \Big[\parci{}{x^\ell} - \G_{\ell \ii m}^{r} E_{\ii c}^{\ii m}  \parci{}{E_{\ii c}^r}\Big]\parci{\Psi}{E_{\ii b}^{k}}.\nonumber
 \end{align*}
It remains to write $\mathscr D \Psi$ and $\mathscr N \Psi$ in terms of the coordinates $(x^i, \tau, e^i_a,e^{\p 0}_a, e^i_{\p 0},e^{\p 0}_{\p 0})$ of $\mathscr G$. Note that to get the formulas \eqref{def_D} and \eqref{N_ex} in the statement of the theorem, it is enough to show that
\bec\label{eq_ts1}
e_{\ii a}^{\ii i} \circ \Phi = E_{\ii a}^{\ii i},
\eec
and  
\bec\label{eq_ts2}
\parcial{(\psi \circ \Phi)}{E^{\ii i}_{\ii a}} = \parcial{\psi}{e^{\ii i}_{\ii a}} \circ \Phi \quad \text{and} \quad  \parcial{(\psi \circ \Phi)}{x^{\ii i}} = \parcial{\psi}{x^{\ii i}} \circ \Phi.
\eec
If $\{e_{\ii a}\}$ is the canonical basis of $\R^{n+1}$, from equation \eqref{uea} we get 
\begin{equation*}
\Phi(\mathscr U) e_{\ii a} = \varpi(\mathscr U E_{\ii a}) = \varpi\Big(E_{\ii a}^{\p I}\parci{}{x^{\p I}}\Big)= E_{\ii a}^{\ii i} \parci{}{x^{\ii i}},
\end{equation*}
which implies \eqref{eq_ts1}. To proceed, observe that equation \eqref{eq_ts1} implies that
\begin{equation*}
\Phi_{\ast}(\parci{}{E^{\ii i}_{\ii a}}) = \parci{}{e^{\ii i}_{\ii a}}.
\end{equation*}
Using this, we compute
\begin{equation*}
\parcial{(\psi \circ \Phi)}{E^{\ii i}_{\ii a}}\bigg|_{\p {\mathscr U}}  = \psi_{\ast \p{\Phi(\mathscr U)}}\left(\Phi_{\ast \p{\mathscr U}} \parcial{}{E^{\ii i}_{\ii a}}\right) = \psi_{\ast \p{\Phi(\mathscr U)}}\left( \parcial{}{e^{\ii i}_{\ii a}}\right) = \parcial{\psi}{e^{\ii i}_{\ii a}}\bigg|_{\p{\Phi(\mathscr U)}},
\end{equation*}
which proves the first part of \eqref{eq_ts2}. Arguing similarly, the identity $x^{\ii i}\circ \Phi=x^{\ii i}$ yields the second part of \eqref{eq_ts2}.

Summing up, we reach the conclusion that
\begin{equation*}
\Delta_{\mathscr{OM} } (\Phi^\ast \psi) = \Phi^\ast  (\mc D^{\p N} \psi) \qquad \text{with} \qquad \mc D^{\p N} \psi = (\mc D  + \mathscr N) \psi + \mathscr E_{\p N},
\end{equation*}
where the operators $\mc D$ and $\mathscr N$ are defined in \eqref{def_D} and \eqref{N_ex}, and $\lim_{N\to \infty} \mathscr E_{\p N} = 0$. Using this specific $\mc D^{\p N}$ in Step 2, the rest of the proof of Theorem \ref{lim_sol_fb} follows now exactly as in the proof of Theorem \ref{lim_sol}.
\end{proof}

\subsection{Conclusion of the proof}
The goal of this section is to finish the proof of Theorem \ref{MainT2}.

Recall that for the Ricci flow we work on the ${\mathrm O}_n$-bundle
\begin{equation*}
\mathscr{O}\to M\times I.
\end{equation*}
The canonical horizontal vector fields on this bundle are defined via the space-time connection \eqref{st_conn}. Namely, there are the $n$ spatial horizontal vector fields $H_i(u)= (ue_i)^\ast$ where $\ast$ denotes the horizontal lift of $ue_i\in (T_x M,g_\tau)$ using the metric $g_\tau$, and the time-like horizontal vector field $D_\tau$ which is the horizontal lift of the vector field $\partial_\tau$. In local coordinates (as always we view ${\mathrm O}_n$-bundle $\mathscr{O}$ sitting inside the corresponding ${\mathrm GL}_n$-bundle $\mathscr{F}$) these vector fields are given by
\bec\lb{rfhor}
H_{a}(u) = e^i_a \left( \frac{\partial}{\partial x^i}-\Gamma_{ij}^k e^j_b \frac{\partial}{\partial e^k_b}\right),
\eec
and
\begin{equation*}
D_\tau(u)  = \parci{}{\tau} -  e^k_a R_{k}^i  \parci{}{e^i_a},
\end{equation*}
cf.~\cite{Ham93,HaNa}. The horizontal Laplacian on $\mathscr{O}$ is given by
\bec\lb{rflap}
\Delta_H=\sum_{a=1}^n H_aH_a.
\eec

\begin{proof}[{Proof of Theorem~\ref{MainT2}.}]
Let $\mathscr{O}\to M\times I$ be the ${\mathrm O}_n$-bundle as above, and let $\mathscr{G}\to M\times I$ be the ${\mathrm GL}_{n+1}$-bundle defined in \eqref{def_gln1}. Recall from the paragraph before the statement of Theorem \ref{lim_sol_fb} that we have the inclusion map
\begin{equation*}
{\mathfrak i}:  \mathscr{O} \to  \mathscr{G},\qquad u \mapsto\left[\begin{array}{c|c}
1 & 0   \\ \hline  0 & u
\end{array}\right]\, .
\end{equation*}
If $\mathbb{Q}_u$ is any subsequential limit as in the statement of Theorem \ref{lim_sol_fb} (existence of such a subsequential limit follows similarly as in Section \ref{sec_subseq_lim}), then it solves the martingale problem for the differential operator $\mc D + \mathscr N$, where $\mc D$ and $\mathscr N$ are defined in \eqref{def_D} and \eqref{N_ex}, respectively. The key is now to observe that $\mathscr N$ restricted to $\mathfrak{i}(\mathscr O)$ vanishes identically, i.e.
\bec\lb{obs_key1}
\mathscr N|_{\mathfrak i(\mathscr O)} \equiv 0,
\eec
and that $\mc D$ restricted to $\mathfrak{i}(\mathscr O)$ is exactly the horizontal heat operator for the Ricci flow, i.e.
\bec \lb{obs_key2}
\mc D|_{\mathfrak i(\mathscr O)} = D_\tau + \Delta_{H},
\eec
where the operator on the right hand side is given by \eqref{rfhor}--\eqref{rflap}.

Consider the projection map
\begin{equation*}
\Pi: \mathscr G\to \mathscr F,\quad \left[\begin{array}{c|c}
e^{\p 0}_{\p 0} & e^{\p 0}_{i}   \\ \hline  e^{\p a}_{\p 0} & u
\end{array}\right] \mapsto u
\end{equation*}
from the ${\mathrm GL}_{n+1}$-bundle $\mathscr G$ to the ${\mathrm GL}_{n}$-bundle $\mathscr F$.
Note that $\Pi(\mathfrak{i}({u}))=u$.

{\bf Claim.} $\Pi(U_s)\in \mathscr{O}$ for all $s\in I$ almost surely. 

\noindent To prove the claim, we use the test function
\begin{equation*}
\psi(U_s) = d_{\mathscr G}(\mathfrak{i}(\Pi(U_s)), \mathfrak{i}(\mathscr O))^2.
\end{equation*}
After multiplication with a suitable cut-off function, we may consider $\psi \in C^\infty_c(\mathscr G)$. Note that Theorem \ref{lim_sol_fb} (a) ensures that $\psi(U_{\p 0})=0$.
Denote by $\mathscr{K}\subset\mathscr{G}$ any compact neighborhood of $\mathfrak{i}(\Pi(U_{\p 0}))$.  By the key observation \eqref{obs_key1}--\eqref{obs_key2}, and since the vector fields $D_\tau$ and $H_a$ are tangential to $\mathfrak{i}(\mathscr{O})$, there exists an open neighborhood $\mathscr{V}\subset\mathscr{G}$ of $\mathfrak{i}(\mathscr{O})$, and a constant $C_{\mathscr{K}}<\infty$, such that
\bec\lb{GW_bound}
(\mc D+\mathscr N)\psi \leq C_{\mathscr{K}} \psi \quad \text{on} \quad \mathscr V\cap \mathscr{K}.
\eec
Now take the first exit time from $\mathscr V\cap \mathscr{K}$ as a stopping time, i.e. set
\begin{equation*}
\sigma= \inf\{s > 0 \, |\, U_s \not\in \mathscr V\} \wedge \inf\{s > 0 \, |\, U_s \not\in \mathscr K\}.
\end{equation*}
Using Theorem \ref{lim_sol_fb} and equation \eqref{GW_bound}, a Gronwall-type argument implies that
\begin{equation*}
\psi(U_{s}) = 0\quad \textrm{ for all $s \leq \sigma$ almost surely.}
\end{equation*}
This in turn means that $\Pi(U_s) \in \mathscr O$ for all $s \leq \sigma$ almost surely. Taking an exhaustion by compact sets $\mathscr K$, we conclude that $\Pi(U_s)\in \mathscr O$ for all $s\in I$ almost surely, as claimed.

Summing up, given $u \in \mathscr O$, we have shown that any subsequential limit of the sequence $\hat{\Pi}_\ast\hat{\Phi}_\ast \overline{\mathbb P}_u^{\p N}$, where $\overline{\mathbb P}_u^{\p N}$ is defined in \eqref{Def_QNu}, concentrates on curves in $P(\mathscr O)$ and solves the martingale problem associated to the operator $D_\tau + \Delta_H$ based at $u$. This coincides with the martingale problem corresponding to the stochastic differential equation \eqref{SDE_RF}, and thus, by uniqueness, we get that the limiting measure must be equal to $\mathbb P_u^{\pa}$ and that the subsequential convergence actually entails full convergence.

Finally, since $\Theta = \Pi \circ \Phi$ this implies the assertion of Theorem \ref{MainT2} as stated in the introduction.
\end{proof}

\section{Application: from Ricci bounds to Ricci flow estimates} \lb{Applic}

In this final section, we explain how one can recover the estimates from \cite{HaNa}.

\begin{proof}[{Proof of Corollary \ref{Cor3}}]
We work with test functions $F:P(\mathscr{M})\to \mathbb{R}$ of the form
\begin{equation*}
F(\mathscr{X})=f(X_{s_1},\ldots,X_{s_k}),
\end{equation*}
where $f:M^k\to \mathbb{R}$ is smooth with compact support, $\frak s=\{0 < s_1< \cdots < s_k < T\}$ are given evaluation times, and $X_{s_i}$ denotes the $M$-component of $\mathscr{X}_{s_i}$.

More systematically, in terms of the commutative diagram
\bec \lb{diag_com1}
\begin{CD}  P(\mathscr M) @>{\mathscr E}_{\frak s}>> {\mathscr M}^k \\  @VV\hat{\pi}V @VV\pi_kV\\ P(\mathbb M) @>\eps_{\frak s}>>  M^k @>f>>  \R,\end{CD}
\eec
where $\mathscr  E_{\frak s}$ and $\eps_{\frak s}$ are the evaluation maps at the $k$ times given by $\frak s$, this means that
\begin{equation*}
F  = (\eps_{\frak s}\circ \hat \pi)^\ast f = (\pi_k \circ \mathscr E_{\frak s})^\ast f .
\end{equation*}
Applying the gradient estimate  \eqref{grad_est_Ric_inf} with $\kappa = C/N$ we get\footnote{To be precise, we use the variant of the gradient estimate for manifolds with boundary established by Wang-Wu \cite{WW16}. Since $\tau(X_s)=s$ almost surely in the limit, the boundary term, which is proportional to the local time spent at the boundary, disappears in the limit.}
\begin{equation} \lb{GS_start}
\limsup_{N\to \infty}\bigg|\nabla_p \int_{P(\mathscr M)} \hat{\pi}^\ast \eps_{\frak s} ^\ast f \, d \mathbb P_p^{\p N}\bigg| \leq \limsup_{N\to \infty} \int_{P(\mathscr M)} \!\! \big|^{\p N}\!\nabla_{0}^{||} ( \mathscr E_{\frak s}^\ast \pi_k^\ast f )\big| \,d  \mathbb P_p^{\p N},
\end{equation}
where we used that $F$ is $\Sigma_T$-measurable to control the error term. Note also that, since our test function doesn't depend on the $\mathbb S^{\p N}$-factor, in \eqref{GS_start} we can replace $\mathbb P_p^{\p N}$ by the averaged probablity measure $\overline{\mathbb{P}}_{{\p (x,T)}}^{\p N}$, where $(x,T)=\pi(p)$.
Using Theorem \ref{MainT1} we compute
\begin{equation} \lb{lim_phi}
\lim_{N\to \infty}\int_{P(\mathscr M)} \hat{\pi}^\ast \eps_{\frak s} ^\ast f \, d \overline{\mathbb{P}}_{{\p (x,T)}}^{\p N}= \int_{P^{\pa}_{\p{T}}(\mathbb M)}  \eps_{\frak s} ^\ast f \, d \mathbb P^{\pa}_{\p{(x, T)}}.
\end{equation}
Combining \eqref{GS_start} and \eqref{lim_phi}, and using lower-semicontinuity, we infer that
\begin{equation} \lb{GS_middle}
\bigg|\nabla_x \int_{P^{\pa}_{\p{T}}(\mathbb M)}  \eps_{\frak s} ^\ast f \, d \mathbb P^{\pa}_{\p{(x, T)}}\bigg| \leq \limsup_{N\to \infty} \int_{P(\mathscr M)} \!\! \big|^{\p N}\!\nabla_{0}^{||} ( \mathscr E_{\frak s}^\ast \pi_k^\ast f )\big| \,d  \mathbb P_p^{\p N},
\end{equation}
where the left hand side is understood as local Lipschitz slope.

To compute the right hand side of \eqref{GS_middle}, it is best to lift things to the frame bundle using the map $\hat{\mathscr P}: P(\mathscr{OM})\to P(\mathscr M)$. By definition of the parallel gradient \cite[Section 6]{Na13} we compute that
\begin{equation*}
\hat{\mathscr P}^\ast \big|^{\p N}\!\nabla_{0}^{||} ( \mathscr E_{\frak s}^\ast \pi_k^\ast f )\big|^2=E_{\frak s}^\ast \sum_{A=0}^{D-1}\bigg|\sum_{i=1}^k H_A^{(i)}(\widetilde{\pi_k^\ast f}) \bigg|^2,
\end{equation*}
where $\widetilde{\pi_k^\ast f}:(\mathscr{OM})^k\to \mathbb{R}$ denotes the invariant lift of $\pi_k^\ast f$, and $E_{\frak s}: P(\mathscr{OM})\to (\mathscr{OM})^k$ is the evaluation map on the frame bundle. Using this and the relation \eqref{rel_P_Pu} we infer that
\begin{equation*}
\int_{P(\mathscr M)} \!\! \big|^{\p N}\!\nabla_{0}^{||} ( \mathscr E_{\frak s}^\ast \pi_k^\ast f )\big| \,d  \mathbb P_p^{\p N}=\int_{P(\mathscr {OM})} \left( E_{\frak s}^\ast\sum_{A=0}^{D-1}\bigg|\sum_{i=1}^k H_A^{(i)}(\widetilde{\pi_k^\ast f}) \bigg|^2\right)^{1/2} d  \mathbb P_{\mathscr U}^{\p N}.
\end{equation*}
Repeating the same argument for the Ricci flow, and using \eqref{rel_P_Pu_par}, we obtain
\begin{equation*}
 \int_{P_{\p{T}}^{\pa}(\mathbb{M})} \!\! \big|^{\pa}\nabla_{0}^{||} (\eps_{\frak s}^\ast f) \big| \, d \mathbb P_{\p{(x, T)}}^{\pa}=\int_{P(\mathscr {O})} \left(e_{\frak s}^\ast\sum_{a=1}^{n}\bigg|\sum_{i=1}^k H_a^{(i)}(\widetilde{ p_k^\ast f}) \bigg|^2\right)^{1/2} d  \mathbb P_{u}^{\pa},
\end{equation*}
where $\widetilde{p_k^\ast f}:\mathscr{O}^k\to\mathbb{R}$ denotes the invariant lift of $f \circ p_k :\mathbb{M}^k\to M^k\to\mathbb{R}$, and $e_{\frak s}:P(\mathscr{O})\to \mathscr{O}^k$ is the evaluation map on space-time. Finally, considering the augmented diagram
\bec \lb{diag_com}
\begin{CD} P(\mathscr{OM}) @>\hat{\mathscr P}>> P(\mathscr M) @>{\mathscr E}_{\frak s}>> {\mathscr M}^k \\ @VV\hat{\Theta}V @VV\hat{\pi}V @VV\pi_kV\\ P(\mathscr F) @>\hat{\frak p}>> P(\mathbb M) @>\eps_{\frak s}>>  M^k @>f>>  \R ,\end{CD}
\eec
and using in particular Theorem \ref{MainT2}, it follows that
\begin{multline*}
\lim_{N\to \infty}\int_{P(\mathscr {OM})} \left(E_{\frak s}^\ast\sum_{A=0}^{D-1}\bigg|\sum_{i=1}^k H_A^{(i)}(\widetilde{\pi_k^\ast f}) \bigg|^2\right)^{1/2} d  \mathbb P_{\mathscr U}^{\p N}\\
=\int_{P(\mathscr {O})} \left(e_{\frak s}^\ast\sum_{a=1}^{n}\bigg|\sum_{i=1}^k H_a^{(i)}(\widetilde{p_k^\ast f}) \bigg|^2\right)^{1/2} d  \mathbb P_{u}^{\pa}.
\end{multline*}
Putting everything together, we conclude that for any cylinder function $\eps_{\frak s} ^\ast f$ we have the estimate
\begin{equation*} 
\bigg|\nabla_x \int_{P^{\pa}_{\p{T}}(\mathbb M)}  \eps_{\frak s} ^\ast f \, d \mathbb P^{\pa}_{\p{(x, T)}}\bigg| \leq  \int_{P_{\p{T}}^{\pa}(\mathbb{M})} \!\! \big|^{\pa}\nabla_{0}^{||} (\eps_{\frak s}^\ast f) \big| \, d \mathbb P_{\p{(x, T)}}^{\pa}.
\end{equation*}
Since cylinder functions are dense in $L^2(P^{\pa}_{\p T}(\mathbb M))$, this proves the corollary.
\end{proof}

\begin{rem}
Similarly, the other main estimates for the Ricci flow from \cite{HaNa}, namely, the quadratic variation estimate
\begin{equation*}
\bigg|\int_{P_{\p{ T}}^{\pa}(\mathbb M)}  \deri{[F_s^{\pa}, F_s^{\pa}]}{s} \, d \mathbb P_{\p{ T}}^{\pa}\bigg| \leq 2 \int_{P_{\p{ T}}^{\pa}(\mathbb M)}  \big|^{\pa}\nabla_{s}^{||}  F\big|^2 \, d \mathbb P_{\p{(x, T)}}^{\pa},
\end{equation*}
the log-Sobolev inequality on path space
\begin{multline*}
\int_{P_{\p T}^{\pa}(\mathbb M)} \Big((F^2)_{\p{s_1}}^{\pa} \log(F^2)_{\p{s_1}}^{\pa} - (F^2)_{\p{s_0}}^{\pa} \log (F^2)_{\p{s_0}}^{\pa}\Big) \, d \mathbb P_{\p{(x, T)}}^{\pa} \\
\leq  4 \int_{P_{\p T}^{\pa}(\mathbb M)} \int_{s_0}^{s_1}  \big|^{\pa}\nabla_{s}^{||}  F\big|^2 \, ds \, d \mathbb P_{\p{(x, T)}}^{\pa}  ,
\end{multline*}
and the sharp spectral gap estimate on path space
\begin{equation*}
\int_{P_{\p T}^{\pa}(\mathbb M)} \big(F_{\p{s_1}}^{\pa} - F_{\p{s_0}}^{\pa} \big)^2 \, d \mathbb P_{\p{(x, T)}}^{\pa} \leq  2 \int_{P_{\p T}^{\pa}(\mathbb M)} \int_{s_0}^{s_1}  \big|^{\pa}\nabla_{s}^{||}  F\big|^2 \, ds \, d \mathbb P_{\p{(x, T)}}^{\pa} 
\end{equation*}
can also be recovered by applying the estimates from \cite{Na13} with $\kappa=C/N$ on Perelman's manifold $\mathscr{M}$ and taking the limit $N\to\infty$ via Theorem \ref{MainT1} and \ref{MainT2}. The computations are almost the same as in the proof of Corollary \ref{Cor3}. Essentially, the only new step is to observe that for $\hat{\pi}^\ast F: P(\mathscr{M})\to P_{\p T}^{\pa}(\mathbb M)\to \mathbb{R}$ the induced martingale
\begin{equation*}
(\hat \pi^\ast F)_s(\mathscr{X}) =\int_{P_{\mathscr{X}_s}(\mathscr M)} \!\!(\hat \pi^\ast F)(\mathscr{X}|_{[0, s]}\ast \tilde{\mathscr{X}}) \, d\mathbb P^{\p N}_{\mathscr{X}_s}(\tilde{\mathscr{X}}) 
\end{equation*}
converges for $N\to \infty$ to $\hat{\pi}^\ast F_s^{\pa}$, where
\begin{equation*}
F_s^{\pa}(\gamma)= \int_{P^{\pa}_{\p{ T-s}}(\mathbb M)} F(\gamma|_{[0, s]} \ast \tilde \gamma)\, d\mathbb P^{\pa}_{\p{\gamma_s}}(\tilde \gamma).
\end{equation*}
\end{rem}


\begin{thebibliography}{99}
\bibitem{Per1} G.~Perelman, The entropy formula for the Ricci flow and its geometric applications, arxiv.org/abs/math.DG/0303109.
\end{thebibliography}
\end{document}